\documentclass{amsart}
\usepackage{amssymb,mathrsfs}%MnSymbol}%,skak}
\usepackage{esint} 
\usepackage[english]{babel}
\def\bN {\mathbf{N}}
\def\bP {\mathbf{P}}

\def\bR {\mathbf{R}}

\def\bT {\mathbf{T}}

\def\bZ {\mathbf{Z}}

\def\R {\mathbf{R}}
\def\T {\mathbf{T}}

\def\cD {\mathcal{D}}

\def\cF {\mathcal{F}}

\def\eps {{\varepsilon}}

\def\d {{\partial}}

\newcommand{\fvp}{f^{\text{VP}}}
\newcommand{\Xvp}{X^{\text{VP}}}
\newcommand{\XiVp}{\Xi^{\text{VP}}}

\newcommand{\rhoVp}{\rho^{\text{VP}}}
\newcommand{\jvp}{j^{\text{VP}}}
\newcommand{\phivp }{\phi^{\text{VP}}}

\newcommand{\fvm}{f^{\text{VM}}}
\newcommand{\Xvm}{X^{\text{VM}}}
\newcommand{\XiVm}{\Xi^{\text{VM}}}

\newcommand{\rhovm}{\rho^{\text{VM}}}
\newcommand{\jvm}{j^{\text{VM}}}
\newcommand{\phivm}{\phi^{\text{VM}}}

\newcommand{\rhoteta}{\rho_{\theta}}
\newcommand{\xiteta}{\xi_{\theta}}

\newcommand{\bdelta}{B_{\delta}}
\newcommand{\bdeltazero}{B^{\eta}_{\delta_0}}
\newcommand{\ba}{\begin{aligned}}
\newcommand{\ea}{\end{aligned}}

\newcommand{\be}{\begin{equation}}
\newcommand{\ee}{\end{equation}}

\newtheorem{Thm}{Theorem}[section]
\newtheorem{Rmk}[Thm]{Remark}

\newtheorem{Lem}[Thm]{Lemma}
\newtheorem{Def}[Thm]{Definition}

\numberwithin{equation}{section}

\begin{document}

\title[Non-relativistic limit with uniform macroscopic bounds]{The non-relativistic limit of the Vlasov-Maxwell system with  uniform macroscopic bounds}

\author[N. Brigouleix]{Nicolas Brigouleix}
\address[N.B.]{CMLS, CNRS, \'Ecole polytechnique, Institut Polytechnique de Paris, 91128
Palaiseau Cedex, France}
\email{nicolas.brigouleix@polytechnique.edu}

\author[D. Han-Kwan]{Daniel Han-Kwan}
\address[D.H.-K.]{CMLS, CNRS, \'Ecole polytechnique, Institut Polytechnique de Paris, 91128
Palaiseau Cedex, France}
\email{daniel.han-kwan@polytechnique.edu}

\begin{abstract}
We study in this paper the non-relativistic limit from Vlasov-Maxwell to Vlasov-Poisson, which corresponds to the regime where the speed of light is large compared to the typical velocities of particles. In contrast with
\cite{Asano-Ukai-86-SMA}, \cite{Degond-86-MMAS}, \cite{Schaeffer-86-CMP} which handle the case of classical solutions, we consider measure-valued solutions, whose moments and electromagnetic fields are assumed to satisfy some uniform bounds. To this end, we use a functional inspired by the one introduced by Loeper in his proof of   uniqueness  for the Vlasov-Poisson system \cite{Loeper-2006}. 
We also build a special class of measure-valued solutions, that enjoy no higher regularity with respect to the momentum variable, but whose moments and electromagnetic fields satisfy all required conditions to enter our framework.
\end{abstract}

 
\maketitle

 \tableofcontents
 
\section{Introduction}
\subsection{The Vlasov-Poisson and Vlasov-Maxwell systems}

This paper is concerned with the non-relativistic limit of the relativistic Vlasov-Maxwell system  towards the classical Vlasov-Poisson system.  These equations govern the evolution of a distribution function $f(t, x, \xi)$ describing a system of charged particles interacting through electromagnetic forces, $t\in [0, \infty)$ being the time variable, $x\in \R^3$ or $\bT^3$ (the $3$-dimensional torus $\bR^3 /\bZ ^3$, equipped with the normalized Lebesgue measure) the space variable, and $\xi\in \bR^3$ the momentum variable. Namely, for any $t\in [0, \infty)$, $f(t,\cdot,\cdot)$ stands for the probability density of particles with phase-space coordinates $(x, \xi)$.%\in\bT_x^3\times \bR_{\xi}^3$.\\

Such systems come from the study of magnetized collisionless plasma. The difference between them lies in the way the electromagnetic force is defined. One can refer to the reference monograph  \cite{Glassey}  of Glassey (in particular Chapters IV, V, VI) for an overview on the background of these kinetic equations.

For the sake of simplicity we consider a system with a single species of particles with charge and mass equal to one, say electrons. For the case of periodic boundary conditions, we shall assume the presence of a background of fixed particles of opposite charge and unit density (typically ions whose mass is much larger than that of electrons). We will denote by $c$ the speed of light  and  by $|v_p|$ the typical velocity of particles, and we will consider
$$
\eps:= \frac{|v_p|}{c} \in (0,1]
$$
as the small parameter converging to $0$, characterizing the so-called non-relativistic limit.

For a given momentum $\xi \in \bR^3$, the particles have in the non-relativistic framework a velocity
$$
v(\xi)=\xi,
$$
while in the relativistic framework a velocity 
\[
v(\xi)=  \frac{\xi}{(1+\eps^2|\xi|^2)^{1/2}}.
\]
Their motion is governed by the Vlasov equation:
\[\d_t f + v(\xi)\cdot \nabla_x f + F\cdot \nabla_{\xi} f = 0,\]
where $F$ stands for the electromagnetic force.

To avoid any confusion the quantity $v(\xi)$ will always designate the relativstic velocity in the following.

Such Vlasov equations are well understood when the self-induced force $F$ enjoys a gain of regularity compared to the distribution function $f$, as in the Vlasov-Poisson system where one only considers the action of the electric field stemming from the Coulomb potential. On the whole space, the existence and uniqueness of global in time classical solutions has been known (either for smooth initial data with compact support or with bounded high order  velocity moments) since the seminal work of Lions and Perthame \cite{LiPer} and Pfaffelmoser \cite{Pfa_JDE_92} (see also \cite{Schaeffer-CPDE-91} for a simplified proof). For periodic boundary conditions, adapting the proof of \cite{Pfa_JDE_92}, the global existence of classical solutions has been established by Batt and Rein  in \cite{Batt-Rein-1991}. 

In our framework, that is 
\begin{itemize}
\item either with no other charged particles in the case of the space domain $\Omega= \R^3$,
\item or
with %periodic boundary conditions and 
a background of massive fixed ions in the case of the space domain $\Omega = \T^3$, \end{itemize}
the Vlasov-Poisson system is the following:
\be\label{vp} \tag{VP}
\left\{
\ba
{}& \d_t \fvp + \xi\cdot \nabla_x \fvp + E(t,x) \cdot \nabla_{\xi}\fvp=0, \quad \,t\in [0, \infty),\, x,\xi\in \Omega \times\bR^3\,,
\\
& E = -\nabla \phivp, \\
&-\Delta \phivp =  \rhoVp-\lambda_\Omega,
\ea
\right.
\ee
where $\rhoVp(t,x)$ stands for the macroscopic density function of particles in the plasma and $ \jvp(t,x)$ for the  current density vector which are defined in the following way:
\[
\rhoVp(t, x)= \displaystyle\int_{\bR^3_\xi} \fvp(t, x, d \xi),\quad
\jvp(t, x) = \displaystyle\int_{\bR^3_\xi} \xi \fvp(t, x, d\xi).
\]
To distinguish between $\bR^3$ and $\bT^3$, we use the notation 
$$\lambda_{\R^3} = 0, \quad \lambda_{\T^3}=1.$$
The system is endowed with an initial condition
$$
f|_{t=0} = f^0,
$$
normalized so that $\int_{\Omega \times \R^3_\xi} f^0 ( dx, d \xi)= 1$.
We will denote  by $\fvp$ a solution to~\eqref{vp}.
\\

When the force $F$ is the Lorentz force produced by an electromagnetic field $(E(t,x), B(t,x))$ whose evolution is governed by the Maxwell equations,  we obtain the Vlasov-Maxwell system:
\be\label{vm} \tag{VM}
\left\{
\ba
{}& \d_t \fvm_\eps + v(\xi)\cdot \nabla_x \fvm_\eps + (E_\eps +\eps v(\xi)\times B_\eps)\cdot \nabla_{\xi}\fvm_\eps=0  \,, 
\\
&\eps \d_t E_\eps = \nabla \times B_\eps - \jvm_\eps, \qquad  \qquad  \,t\in [0, \infty),\, x, \xi\in \Omega \times\bR^3\,,
\\  
&\nabla \cdot E_\eps = \rhovm_\eps - \lambda_\Omega,
\\
&\eps\d_t B_\eps = - \nabla \times E_\eps ,\\ 
&\nabla \cdot B_\eps =0, 
\ea
\right. 
\ee
where this time
\[
\rhovm_\eps(t, x)= \displaystyle\int_{\bR^3_\xi} \fvm_\eps(t, x, d \xi),\quad
\jvm_\eps(t, x) = \displaystyle\int_{\bR^3_\xi} v(\xi) \fvm_\eps(t, x, d\xi).
\]
The system is endowed with an initial condition
$$
f^\eps|_{t=0} = f^0, \quad  E_{\eps}|_{t=0} = E^0_\eps, \quad B_{\eps}|_{t=0}= B^0_\eps.
$$
We will denote by $(f_{\eps}^{\text{VM}}, E_{\eps}, B_{\eps})$ a  solution to~\eqref{vm}. We will  often omit to index the quantities involved in~\eqref{vm} by $\eps$ when not needed, in order to enlighten the notations.
 
In the following we will use the potential formulation of the Maxwell equations, sometimes called the Lienard-Wiechert formulation (see for instance Section $2$ of \cite{Bouchut-Golse-Pallard-03ARMA} for a proof of the equivalence of the two formulations, and  \cite{LL}  for the physical point of view) that consists in introducing the potentials $(\phi_\eps, A_\eps)$ solving
\be\label{potential} 
\left\{
\ba
{}& -\Delta  \phivm_\eps = \rhovm_\eps -\lambda_\Omega ,\\
&\eps^2\partial_t^2A_\eps-\Delta A_\eps=\eps \bP(\jvm_\eps),
\ea
\right.
\ee
where $\bP$ stands for the Leray projection, that is the projection from the set of vector fields whose components belong to $L^2(\Omega)$ onto the subspace of divergence free vector fields. It appears in the equations because to have a unique pair of potentials $(\phi_\eps, A_\eps)$ corresponding to the electromagnetic field $(E_\eps,B_\eps)$ we need an additional gauge condition. We make the choice of the Coulomb gauge condition
$$
\nabla\cdot A_\eps =0,
$$
which is encoded in the equations thanks to the operator $\bP$. Then one obtains $E_\eps$ by the formula
\begin{equation}
\label{formula-E}
 E_\eps = -\nabla \phivm_\eps - \eps \d_t A_\eps,
\end{equation}
while $B_\eps$ is recovered 
\begin{itemize}
\item for $\Omega = \R^3$, thanks to the equation
\begin{equation}
\label{formula-B-R}
B_\eps = \nabla\times A_\eps,
\end{equation}
\item for $\Omega = \T^3$, this requires to take into account the mean space  value, which results in the equation
\begin{equation}
\label{formula-B-T}
B_\eps = \nabla\times A_\eps + \int_{\bT^3} B_\eps  \, dx.
\end{equation}
\end{itemize}
We shall also endow~\eqref{potential} with initial conditions $(A_\eps |_{t=0}, \eps\partial_t A_\eps |_{t=0})$ that must be compatible with $E_{\eps}|_{t=0} = E^0_\eps,  B_{\eps}|_{t=0}= B^0_\eps$.

In this paper, as we shall manipulate solutions to the Vlasov-Maxwell and Vlasov-Poisson systems, we will index the quantities $f$, $\rho$, $j$, $\phi$ etc. by $\text{VP}$ or $\text{VM}$ when needed to make the distinction.

\subsection{The non-relativistic limit}
Formally taking $\eps=0$ in the Vlasov-Maxwell system \eqref{vm}, one almost readily obtains the Vlasov-Poisson system \eqref{vp}.
The general goal is to determine, given a sequence $(\fvm_\eps)_\eps$ of weak solutions to \eqref{vm}, whether it converges (in a sense to be made precise) to a  solution $\fvp$ to \eqref{vp}.
This is what we refer to as the non-relativistic limit.

Until now this problem has been tackled in the pioneering (simultaneous and independent) works of Asano-Ukai \cite{Asano-Ukai-86-SMA}, Degond \cite{Degond-86-MMAS} and Schaeffer \cite{Schaeffer-86-CMP}, which all concern classical solutions to \eqref{vm} and \eqref{vp} on $\Omega= \bR^3$. They rely on high order (Sobolev or Lipschitz) estimates  on the difference $\fvp - \fvm_{\eps}$, i.e. between two classical solutions stemming from the same smooth initial distribution function. These non-relativistic limits  require at least Lipschitz uniform regularity and boundedness for the distribution function (with respect to $x$ and $\xi$) and for the electromagnetic field (with respect to $x$). We would like to be able to relax such an assumption. Let us also mention that the non-relativistic limit  has also been treated for the Vlasov-Nordtr\"om system in \cite{Calo-Lee-2004}, for the Vlasov-Maxwell system in lower dimensions in \cite{Lee-2004} and more recently in \cite{Schaeffer-Wu}. Large time (with respect to $1/\eps$) estimates are studied in \cite{Pallard-08-AA} for  $\Omega= \bR^3$ and in \cite{HKN}, \cite{HanKwan-N.-R.2018CMP}, in relation with  Penrose stability issues for 
$\Omega= \bT^3$.

We observe that the structure of the  proofs in \cite{Asano-Ukai-86-SMA}, \cite{Degond-86-MMAS}, \cite{Schaeffer-86-CMP} is in fact somehow similar to the one used by Robert in \cite{Robert-97} to prove a uniqueness result for the Vlasov-Poisson system, when the initial distribution function $f^0$ is a bounded measurable and compactly supported function. More recently, this assumption on the initial data has been weakened by Loeper in \cite{Loeper-2006} who was able to handle measure-valued solutions without compact support. 
Using tools from optimal transport, he proved,  given an initial condition, the uniqueness of weak solutions to the Vlasov-Poisson system which have a bounded macroscopic density. Let us mention that  this result  has been refined by Miot  \cite{Miot-2016} and Holding-Miot \cite{Holding-Miot-2018}.

Our idea in this work is to adopt the same point of view as Loeper in \cite{Loeper-2006}, that is to consider measure-valued solutions to~\eqref{vm} and \eqref{vp} and prove the non-relativistic with conditional bounds on macroscopic quantities (moments or force fields). Loosely speaking, we shall prove that given a fixed initial distribution function $f^0$ and
\begin{itemize}

\item a sequence of solutions $(\fvm_\eps)$  to~\eqref{vm} with a bounded density $\rhovm$ that is uniformly bounded in $\eps$ and with a higher order moment and and electromagnetic field of controlled growth in $\eps$,

\item a solution $\fvp$ to~\eqref{vm} with a uniformly bounded density $\rhoVp$,

\end{itemize} 
together with some other milder conditions to be later stated,
the sequence $(\fvm_\eps)$ must converge to $\fvp$ in the weak-$\star$ sense of measures (which will be quantified using a Wasserstein distance).
Perhaps surprisingly, this framework allows to consider solutions to the Vlasov-Maxwell system with electromagnetic fields that are not even bounded in $\eps$.

To do so, we will derive an Osgood-type inequality on a functional involving the two Lagrangian flows associated to the solutions $\fvp$ and $\fvm_\eps$, an argument related to the proof developed by Dobrushin in \cite{Dobrushin-1979} for the purpose of the Mean-Field limit and then used by Loeper in \cite{Loeper-2006} to prove his uniqueness result for the Vlasov-Poisson system.
Finally, we mention that the idea of using Wasserstein stability estimates for studying singular limits of Vlasov equations has also recently been used in \cite{Han-Kwan-Iacobelli-CMS}, \cite{Han-KwanIcaobelli-JDE} (in the context of the so-called quasineutral limit).

\subsection{Main results}
Before being allowed to state our main result,  we have to specify the notions of weak solutions we will manipulate.

We denote by $\mathcal{P}_2(\Omega\times\bR^3)$ the set of probability measures on $\Omega \times\bR^3$, with finite first two moments, endowed with the standard weak-$*$ topology. We look for solutions to \eqref{vp} and \eqref{vm} such that
$$    f|_{t=0}= f^0.
$$
with $f^0 \in \mathcal{P}_2(\Omega\times\bR^3)$.
\begin{Def}\label{defsolweakvp}(Weak solution to (VP)). For $T>0$, we will call $\fvp$ a weak solution to \eqref{vp} on $[0,T)$ associated to the initial condition $f^0$ if
\begin{itemize}
    \item $f\in C([0,T), \mathcal{P}_2( \Omega \times\bR^3)-w*)$,
    \item $-\nabla \phivp \in L^1(0,T; C(\Omega))$,
    \item for all test functions $\varphi \in C_c^{\infty}([0,T)\times \Omega\times\bR^3)$,
\begin{equation}
    \int_{[0,T)\times\Omega\times\bR^3}\left(\d_t\varphi + v(\xi)\cdot\nabla_x\varphi- \nabla_x\phivp\cdot\nabla_{\xi}\varphi\right)  \, f(t, d\xi, dx)  \, dt = - \int_{ \Omega \times\bR^3}  \varphi|_{t=0} \, f^0(d\xi, dx),
\end{equation}
\item for all $t\in[0,T)$, $\phivp(t)$ solves the Poisson equation:
$$
- \Delta \phivp = \rhoVp -\lambda_\Omega.
$$
\item  for all $t\in[0,T)$, the density of current $\jvp$ satisfies the conservation law:
$$
 \d_t \jvp + \nabla_x : \langle f , |\xi|^2 \rangle = - \nabla_x \phivp \rhoVp .
$$
\end{itemize}
\end{Def}

\begin{Rmk}
We should mention that in the following we only consider solutions such that the macroscopic density $\rhoVp$ belongs to $L^\infty(0,T; L^1 \cap L^{\infty}(\Omega))$ and therefore $\nabla \phivp$ will actually be log-Lipschitz.
\end{Rmk}

\begin{Def}\label{defsolweakvm} (Weak solution to (VM)). For $T>0$, we will call $(\fvm,E,B)$ a weak solution to \eqref{vm} on $[0,T)$ associated to the initial condition $(f^0,E^0,B^0)$ if %in $\mathcal{D}'([0,T)\times \Omega\times\bR^3)$,  if:
\begin{itemize}
    \item $\fvm\in C([0,T); \mathcal{P}_2(\Omega\times\bR^3)-w*)$,
    \item  $E,B \in L^1(0,T; C^0\cap L^\infty(\Omega)) \cap  L^1(0,T; BV_{\text{loc}}(\Omega))$,
    \item for all test functions $\varphi \in C_c^{\infty}([0,T)\times\Omega\times\bR^3)$,
\begin{align*}
    &\int_{[0,T)\times\Omega\times\bR^3}\left(\d_t\varphi + v(\xi)\cdot\nabla_x\varphi +  E\cdot \nabla_{\xi}\varphi  + \eps v(\xi) \times B \cdot \nabla_\xi \varphi \right) \, f(t, d\xi, dx)  \, dt \\
&\qquad \qquad \qquad  = - \int_{\Omega\times\bR^3}  \varphi|_{t=0}\, f^0( d\xi, dx),
\end{align*}
\item for all test functions  $\psi \in C_c^{\infty}([0,T)\times\Omega,\R)$, $\Psi \in C_c^{\infty}([0,T)\times\Omega,\R^3)^3$, %for all $t\in[0,T)$, $(E,B)$ solve the Maxwell equations
$$
\ba
& \eps \langle  E, \d_t \Psi \rangle_{\mathcal{D}',\mathcal{D}}+ \langle   \nabla \times B, \Psi \rangle_{\mathcal{D}',\mathcal{D}}   - \langle \jvm, \Psi \rangle_{\mathcal{D}',\mathcal{D}} =  - \eps \langle E^0, \Psi|_{t=0}\rangle_{\mathcal{D}',\mathcal{D}} \,,
\\  
&\langle \nabla \cdot E,\psi \rangle_{\mathcal{D}',\mathcal{D}} = \langle\rhovm - \lambda_\Omega, \psi\rangle_{\mathcal{D}',\mathcal{D}}
\\
&\eps \langle \d_t B, \Psi \rangle_{\mathcal{D}',\mathcal{D}}  -\langle \nabla \times E, \Psi \rangle_{\mathcal{D}',\mathcal{D}} = \eps  \langle B^0, \Psi|_{t=0}\rangle_{\mathcal{D}',\mathcal{D}}\\ 
&\langle \nabla \cdot B, \psi \rangle_{\mathcal{D}',\mathcal{D}} =0.
\ea
$$
\end{itemize}
\end{Def}

\begin{Rmk}
Let us  provide some remarks about Definition~\ref{defsolweakvm}.

\begin{itemize}
 
 \item Again, we shall only consider solutions such that the density $\rhovm$  belongs to $L^\infty(0,T; L^1 \cap L^{\infty}(\Omega))$. The assumption on the regularity of $E$ is therefore only a condition on $\eps\d_t A$.% since the assumption $\rhovm\in L^1\cap L^{\infty}([0,T)\times\bT^3)$ implies that $\nabla_x \phivm$ is in $ L^{\infty}([0,T)\times\bT^3)$.

  \item The BV-regularity condition on $E$ and $B$ is \emph{ad hoc} in order to be able to define the Lagrangian flow thanks to the theory developed by Ambrosio \cite{Ambrosio-2004-Inventionnes} (generalizing the theory of DiPerna and Lions \cite{DipernaLions-89-ODEs}). Note however that this regularity is not asked to be uniform with respect to $\eps$.

  \item  To ensure this regularity, we can for instance impose some Sobolev regularity for $\jvm$. Thanks to the wave equation satisfied by $A$ in \eqref{potential} and the Sobolev embedding, it is for instance enough to set the condition $\jvm \in L^{\infty}(0,T;H^{s}(\Omega))$ for some $s>3/2$.

\end{itemize}
\end{Rmk}

Before stating our main result we also need to specify further assumptions that the solutions we are considering need to satisfy.

\begin{Def}
\label{sol-vp}
We say that a weak solution $\fvp$ to the Vlasov-Poisson system in the sense of Definition \ref{defsolweakvp} is suitable if
 the macroscopic density and the fourth moment of $\fvp$ are bounded in the following sense:
$$
\|\rhoVp\|_{L^{\infty}(0,T; L^1\cap L^\infty(\Omega))} < +\infty,  \qquad \left\|\int_{\bR^3}|\xi|^4\fvp(t,x,d \xi) \right\|_{L^{\infty} (0,T;L^1(\Omega))} < +\infty.
$$
%\end{itemize}

\end{Def}

Such a solution is actually unique according to \cite{Loeper-2006} (recalled in Theorem~\ref{thm-loeper} below).

Finally, some normalization conditions are required for the initial conditions that we are going to consider. These are stated in the following definition.

\begin{Def}
\label{cond-initial}
We say that the initial data $(f^0, E^0_\eps, B^0_\eps) \in \mathcal{P}_2 (\Omega\times\bR^3)\times \mathcal{D}'(\Omega)^2$ are normalized if the following conditions hold.
   First of all,
    $$
\nabla \cdot E^0_\eps = \int_{\bR^3} f^0(x,d \xi) -\lambda_\Omega, \quad { and }\quad \nabla\cdot B^0_\eps = 0  \quad \text{ in } \cD'(\Omega).
    $$
In the case $\Omega = \T^3$, we furthermore  ask that 
\begin{itemize}
    \item the spatial mean of the current density is initially assumed to be zero\footnote{Note that this assumption on the current density can be ensured by a Galilean change of frame.}    
    $$
 %   \int_{\bT^3\times\bR^3 }f^0(x,\xi)  \, d\xi dx =1 \quad \text{ and } \quad 
 \int_{\bT^3\times\bR^3}\xi f^0(dx,d \xi) = 0;
    $$
    \item the spatial mean value of $E^0_\eps$ satisfies
     $$
    \langle E^0_\eps \rangle := \int_{\bT^3} E^0_\eps(x) \, dx = 0. 
        $$
\end{itemize} 
    
\end{Def}

We are finally in position to state our main result:
\begin{Thm}\label{theo1}
Let $f^0$ in $\mathcal{P}_2(\Omega\times\bR^3)$, $(E^0_\eps,B^0_\eps) \in \mathcal{D}'(\Omega)^2$ be normalized  initial data in the sense of Definition~\ref{cond-initial}.
Let $T>0$ (independent of the parameter $\eps$) and assume that $\fvp$ and $(\fvm_{\eps},E_{\eps}, B_{\eps})$ are weak solutions to respectively \eqref{vp} and \eqref{vm} on the interval of time $[0,T]$, with respective  initial data $f^0$ and $(f^0, E^0_\eps,B^0_\eps)$. %, with $T$ %, and satisfying respectively the assumptions $H_1$ and $H_2$.\\
Assume furthermore that $\fvp$ is a suitable solution in the sense of Definition~\ref{sol-vp}.

Finally assume that there exists $C_0>0$  
and $(\alpha, \beta, \gamma_1, \gamma_2) \in [0,1)^4$ such that
\begin{equation}
\kappa:= \min(\alpha- (\beta+2 \gamma_2), 1-(\gamma_1+\gamma_2))>0,
\end{equation}
so that
\begin{itemize}
    \item the macroscopic density $\rhovm$ is uniformly bounded in $\eps$:
    \begin{equation}
    \label{bound-rho-vm}
    \|\rhovm_\eps\|_{L^{\infty}(0,T;L^1 \cap L^{\infty}(\Omega))} \leq C_0,
\end{equation}
    \item   the moment of order $\alpha$
    $$m_\alpha (t,x) := \int_{\bR^3}\fvm_\eps(t,x,d \xi) |v(\xi)|^{\alpha},$$ 
    has a uniform bound that has a controlled growth in $\eps$:
    \begin{equation}
    \label{malpha}
  \|m_\alpha \|_{L^{\infty}(0,T;L^{\infty}(\Omega))} \leq C_0\eps^{-\beta},
\end{equation}
    \item the longitudinal electric and magnetic fields have  a $L^2$ norm that has  a controlled growth in $\eps$:
\begin{equation}
\label{growth-fields}
\begin{aligned}
\|\eps \partial_t A_\eps\|_{L^{\infty}(0,T;L^2(\Omega))} &\leq  C_0 \eps^{-\gamma_1},\\
 \|B_\eps\|_{L^{\infty}(0,T;L^2(\Omega))}  &\leq  C_0 \eps^{-\gamma_2}.
\end{aligned}
\end{equation}
\end{itemize}
Then there exist $\eps_0>0$ a constant $C>0$ depending only on the initial data such that\footnote{Here $W_2$ stands for the Wasserstein-2 distance on which we refer to Section~\ref{sec-Wasserstein} for a definition and some properties.} for all $\eps \in (0,\eps_0]$ and for all $t$ in $[0,T]$,
\begin{equation}\label{inegalitenonrelativiste}
W_2(\fvp,\fvm_{\eps})(t) \leq  \left(C\left(1+T\right)^2\eps^\kappa \right)^{\exp\left(-C(1+T)^2\right)}.
\end{equation}

\end{Thm}

The proof of Theorem~\ref{theo1} is based on the study of the functional
\be
\label{Q-intro}
Q(t) : = \frac{1}{2}\int_{\Omega\times\bR^3}f^0(dx,d\xi)\left(|\Xvp-\Xvm|^2 + |\XiVp-\XiVm|^2\right)
\ee
where, loosely speaking, $(\Xvp, \XiVp)$ (resp. $(\Xvm, \XiVm)$) stand for the characteristic curves (more precisely the lagrangian flow) associated to Vlasov-Poisson (resp. Vlasov-Maxwell).

The core of the proof will consist in proving the following Osgood-type inequality. For all $t\in [0,T]$,
\begin{equation}
\label{eq-Osgood-intro}
Q(t)\leq C\left(1+T\right)^2 \eps^{\kappa} + \int_0^t C\left(1+T\right)^2 Q(s)\left(1+ \log^+\left(\frac{1}{Q(s)}\right)  \right) ds,
\end{equation}
from which we can deduce the stability estimate~\eqref{inegalitenonrelativiste}. \\

Of course, it is possible to apply Theorem~\ref{theo1} to the strong solutions built in \cite{Asano-Ukai-86-SMA}, \cite{Degond-86-MMAS}, \cite{Schaeffer-86-CMP}, for which we already know that the non-relativistic limits holds. 
But Theorem~\ref{theo1} is designed to handle other types of solutions.
It is the purpose of a second result (see Theorem~\ref{thm2} later in this paper) to build measure-valued weak solutions to~\eqref{vm} and~\eqref{vp} that do not enjoy higher regularity with respect to $\xi$, but whose moments and electromagnetic field satisfy all required conditions.

\subsection{Comments on Theorem~\ref{theo1}}
A few comments on the statement of Theorem~\ref{theo1} are in order.

\bigskip

    \noindent{\bf 1.}  Contrary to the classical results of \cite{Asano-Ukai-86-SMA}, \cite{Degond-86-MMAS}, \cite{Schaeffer-86-CMP}, we  obtain in the stability estimate~\eqref{inegalitenonrelativiste}  a polynomial rate of convergence  whose exponent decreases exponentially fast to zero with the time running.

It is not clear whether the method of \cite{Holding-Miot-2018} instead of that of \cite{Loeper-2006} can be adapted to study the non-relativistic limit: if it were the case, this would allow to improve the rates of convergence in the stability estimate~\eqref{inegalitenonrelativiste}. This is left to future studies.
       \medskip

           \noindent{\bf 2.}  As already mentioned, Theorem~\ref{theo1} allows electromagnetic fields that may blow up in $\eps$. No well-prepared assumption is either required.
            Note however that the convergence result only concerns distribution functions and not these fields.
       
       \medskip
       
        \noindent{\bf 3.}   One can observe that the condition on the moments of the solution to the Vlasov-Maxwell system $\fvm$ is (perhaps as expected) more restrictive than the one imposed on $\fvp$ (which corresponds to the criterion of Loeper \cite{Loeper-2006} in view of uniqueness), since we require a control of a higher order moment $m_\alpha$. But  the parameter $\alpha
    $ can be taken arbitrarily close to $0$ and some growth in $\eps$ is even permitted. The price to pay is that the rate of convergence in~\eqref{inegalitenonrelativiste} gets  deteriorated at the same time.

\medskip

    \noindent{\bf 4.}  Our framework shows little dependence on the space dimension: all statements could be modified to handle other dimensions than $3$. This aspect differs from the classical approaches on $\Omega=\bR^d$ where explicit formulas for solutions to wave equations are used and are thus dependent on the dimension $d$ (in particular on its evenness or oddness).

\medskip

    \noindent{\bf 5.}    If additionally $f^0$ is in $L^1 \cap L^{\infty}(\Omega \times\bR^3)$ and $E^0_\eps, B^0_\eps \in L^2(\Omega)$, then for all $\eps \in (0,1]$, one can build a global weak solution $\fvm_\eps$ to~\eqref{vm} satisfying an energy inequality (see Theorem~\ref{thm-DPL} below, and Proposition $1.6$ of \cite{bouchut-golse-pulvirenti-2000}). If we additionally  assume that the initial energy is uniformly bounded in $\eps$, that is to say 
    $$
    \ba
 &\int_{\Omega \times \R^3} \frac{1}{\eps^2}\left(\sqrt{1+ \eps^2 |\xi|^2}-1\right) f^0 \, d\xi dx + \frac{1}{2} \int_{\Omega} \left( |E^0_\eps|^2 + |B^0_\eps|^2 \right) \, dx  \\
   &\qquad \leq \int_{\Omega \times \R^3}  |\xi|^2 f^0 \, d\xi dx + \frac{1}{2} \int_{\Omega} \left( |E^0_\eps|^2 + |B^0_\eps|^2 \right) \, dx \leq C_0,
   \ea
    $$
   we obtain the existence of $C_1>0$ such that for all $\eps \in (0,1]$,
    $$
   \| |\xi|\fvm_\eps \|_{L^\infty(0,+\infty; L^1(\Omega \times\bR^3))} + \| E_\eps \|_{L^\infty(0,+\infty; L^2(\Omega))} + \|B_\eps \|_{L^\infty(0,+\infty; L^2(\Omega))} \leq C_1.
    $$
    This means that the control~\eqref{growth-fields} with $\gamma=0$ is automatically ensured for such solutions. 

    \medskip
    
            \noindent{\bf 6.}    In the case $\Omega= \bT^3$, it is also possible to consider  initial data such that that the spatial mean-value of the initial electric field satisfies
        $$
       |\langle E^0_\eps \rangle| \leq C_0 \eps^{\delta},
        $$
    for some $\delta >0$. The rate of convergence in~\eqref{inegalitenonrelativiste} can get worse  if $\delta$ is close to $0$.
        \medskip

        \noindent{\bf 7.}  It is finally possible to consider initial conditions that fully depend on $\eps$, i.e. $(f^0_\eps, E^0_\eps, B^0_\eps)$ for Vlasov-Maxwell and an initial condition $f^0$ for Vlasov-Poisson. Loosely speaking, the final stability estimate~\eqref{inegalitenonrelativiste} is then replaced by 
        \begin{equation}\label{inegalitenonrelativiste-eps}
W_2(\fvp,\fvm_{\eps})(t) \leq C \left(C\left(1+T\right)^2\eps^\kappa  + C W_2(f^0_\eps, f^0)\right)^{\exp\left(-C(1+T)^2\right)},
\end{equation}
assuming $W_2(f^0_\eps, f^0)$ is small enough. To obtain such a result, this requires to modify the functional $Q$ introduced in~\eqref{Q-intro} (see for instance the proof of Theorem 3.1 in \cite{Han-KwanIcaobelli-JDE}).

\subsection{The Cauchy problem for the Vlasov-Maxwell and Vlasov-Poisson systems: a short review}

The study of the Cauchy problem (either for classical or weak solutions) for the Vlasov-Poisson and the Vlasov-Maxwell system has a long history. We will only (quickly) review some aspects that are pertaining to this work.
\\

\noindent {\bf Vlasov-Poisson.} For what concerns classical solutions, we have already discussed the landmark works \cite{LiPer}, \cite{Pfa_JDE_92} and \cite{Schaeffer-CPDE-91}, \cite{Batt-Rein-1991} (see also the earlier important works \cite{UO} and \cite{BD}).
On the other hand, Arsenev in \cite{Arsenev} built the first global weak solutions to~\eqref{vp}. Let us also state  the uniqueness result of \cite{Loeper-2006} to which this work is related.
\begin{Thm}
\label{thm-loeper}
(Loeper \cite{Loeper-2006}) Given $f^0$ in $\mathcal{P}_2 (\bR^3\times\bR^3)$, there exists at most one weak solution $f$ to~\eqref{vp} 
such that
$$
\left\|\int_{\bR^3}f(t,x,d \xi)\right\|_{L^{\infty}([0,T)\times\bR^3)} < + \infty.
$$
\end{Thm}
The propagation of moments is also an important issue that was studied in \cite{LiPer} on $\bR^3$.
Recently, Pallard in \cite{Pallard-2012-CPDE} was able to prove the propagation of moments on $\bT^3$ (see also \cite{Pallard-sima}).
\begin{Thm}(Pallard \cite{Pallard-2012-CPDE})
Given $k>14/3$ and a non-negative initial data $f^0 \in L^1 \cap L^{\infty} (\bT^3 \times \bR^3)$ such that:
$$
\int_{\bT^3 \times \bR^3} |\xi|^k f^0(x,\xi) \, d\xi dx < + \infty,
$$
there exists a weak solution $f\in C([0,T), \mathcal{P}_2(\bT^3\times\bR^3)-w*)$ to the Cauchy problem for the Vlasov-Poisson system \eqref{vp} such that for any $T>0$ we have the following propagation of the moments:
$$
\int_{\bT^3 \times\bR^3} |\xi|^k f(t,x,\xi) \, d\xi dx< + \infty.
$$
\end{Thm}

\bigskip

\noindent {\bf Vlasov-Maxwell.} The theory of (local) classical solutions to the Vlasov-Maxwell system (seen as a quasi-linear equation) was first developed in \cite{Woll-84}, \cite{Asano}, \cite{Glassey-Strauss86ARMA}. In particular in \cite{Glassey-Strauss86ARMA}, Glassey and Strauss found a criterion for the formation of possible singularities: loosely speaking, they can occur only at large velocities. This was later revisited in \cite{Bouchut-Golse-Pallard-03ARMA}, \cite{Kla.-Sta.-02CPAA}.  We also refer to \cite{Pallard-CMS} and \cite{LS} (and references therein) for some recent developments.
Global existence is only known for special cases, in particular for small data (see \cite{Gla.-Scha.-88-CMP}, \cite{Schaeffer-IUMJ}, \cite{Big}).
The existence of global weak solutions was proved in the landmark work by DiPerna and Lions in \cite{DiPerna-Lions-89-CPAM}.
\begin{Thm}\label{thm-DPL}(DiPerna-Lions \cite{DiPerna-Lions-89-CPAM})
Let $ f^0\in L^1\cap L^\infty(\bR^3 \times \bR^3)$ be a non-negative function that satisfies:
$$
\int_{\bR^3 \times \bR^3} |\xi|^2 f^0(x,\xi) \, d\xi dx < + \infty.
$$
Let $E^0, B^0 \in L^2(\bR^3)$ satisfy the following compatibility conditions on the initial data:
$$
\nabla_x\cdot E^0 = \int_{\bR^3} f^0(x,\xi) \, d\xi, \quad { and }\quad \nabla\cdot B^0 = 0  \quad \text{ in } \cD'(\bR^3).
$$
Then there exist $f\in L^{\infty}(0,+\infty; L^1 \cap L^\infty(\bR^3 \times \bR^3))$ and $E,B \in  L^{\infty}(0,+\infty; L^2(\bR^3))$ which satisfy  \eqref{vm} in the sense of distributions.
\end{Thm}
As far as we know, propagation of moments for solutions to~\eqref{vm} remains largely open.

\subsection{Organization of the paper} The paper is organized as follows. In the next Section \ref{sec-Wasserstein} we recall some basic definitions and facts about the Wasserstein distance in order to set up the framework we will work in and also recall a very useful property proved by Loeper in \cite{Loeper-2006}.  Then we proceed to the proof of Theorem~\ref{theo1} in Section \ref{sec-prooftheo1}. Namely we provide an Osgood inequality for the functional $Q$. Finally in Section \ref{sec-multi}  and in Section \ref{sec-prooftheo2}, in the case $\Omega= \bT^3$, we construct a special class of measure-valued solutions, that have no higher regularity in $\xi$ but that are very regular in $x$, namely real-analytic.
This is based on a multifluid representation (introduced by Grenier in \cite{Grenier-96}), with analyticity regularity in space variable but only measure in the momentum variable.\\

Throughout this paper, $C$ will designate a positive constant depending on the initial data  but independent of the parameter $\eps$, that may change from line to line.

 \bigskip
 
 \noindent {\bf Acknowledgements.} Partial support of the grant ANR-19-CE40-0004 is acknowledged.

\section{Definitions, notations and classical results about the Wasserstein distance $W_2$}
\label{sec-Wasserstein}

This section is devoted to the exposition of a few classical results from  Optimal Transportation Theory (see \cite{Villani-2003-TopicsinOpti} for an overview of the tools needed here) and the links between the Wasserstein distances and the $H^{-1}$-Sobolev norm. For the proofs of the following theorems we refer to Section $2$ of \cite{Loeper-2006}.

\begin{Def}
Let $X$ and $Y$ be two polish spaces. Let $\rho_1$, $\rho_2$ be two Borel probability measures on respectively $X$ and $Y$. We define the Wasserstein distance of order $2$ between $\rho_1$ and $\rho_2$, denoted $W_2(\rho_1, \rho_2)$, by:
$$
W_2(\rho_1, \rho_2) = \underset{\gamma}{\inf} \left(\int_{X\times Y}d(x,y)^2 d\gamma(x,y) \right) ^{1/2},
$$
where the $\inf$ runs over the set of probability measures $\gamma$ on $X\times Y$ whose marginals $\mathbb{P}_x \gamma$ and $\mathbb{P}_y \gamma$ are equal respectively to $\rho_1$ and $\rho_2$. 
\end{Def}

\begin{Rmk} Let us state some remarks about the previous definitions.
 \begin{itemize}
\item We do not need this degree of generality for our purpose, in the following $X$ and $Y$ will always be either $\bT^3$ or $\bR^3$.
\item The Wasserstein distance of order $p$ would have been defined in the same way, only replacing $d(x,y)^2$ by $d(x,y)^p$ but we restrict ourselves to the case $p=2$.
\item There is an important relation between this distance and the optimal transportation theory. This is what enables to relate the distance $W_2$ and the $H^{-1}$-norm, a relation described in the next proposition. The proof and the intermediary lemmas that lead to this result for measure on $\bR^3$ are detailed in Section $2$ of \cite{Loeper-2006} (they adapt to $\bT^3$ with minor changes), based on the seminal results from  optimal transportation theory by Brenier \cite{Brenier-1991} and McCann and Gangbo \cite{Gangbo-96}.
\end{itemize}
\end{Rmk}

\begin{Thm} \label{theo3}(Loeper \cite{Loeper-2006}.) 
Let $\rho_1$, $\rho_2$ be two probability measures on  $\Omega$ with $L^{\infty}$ density with respect to the Lebesgue measure. Let $\psi_i$, $i=1,2$, solve:
\begin{align*}
- \Delta \psi_i = \rho_i -1, \quad \text{on } \Omega.
\end{align*}
Then
$$
\|\nabla\psi_1- \nabla \psi_2 \|_{L^2(\Omega)} \leq \left(\max\{\|\rho_1\|_{L^{\infty}}, \|\rho_2\|_{L^{\infty}}\}\right)^{1/2} W_2(\rho_1, \rho_2).
$$
\end{Thm}

\begin{Def}
Let $\rho_1$ be a Borel probability measure on $\Omega \times \bR^3$ and $T: \Omega\times \bR^3 \rightarrow \Omega\times  \bR^3$ be a measurable mapping. The push-forward of $\rho_1$ by $T$ is the measure $\rho_2$ defined by
$$
\forall B \subset \Omega\times \bR^3 \, \text{Borel}, \quad \rho_2(B) = \rho_1(T^{-1}(B)).
$$
We will use the notation $ \rho_2 = T \# \rho_1$.
\end{Def}

We shall now give a useful remark in view of the estimates of section $3$.

\begin{Rmk}
Let $(\Omega_0, \mu)$ be a probability space , and consider $X_1, X_2$ two mappings from $(\Omega_0, \mu)$ to $\Omega \times \bR^3$. If $X_1 \#d\mu = \rho_1$, $X_2 \# d\mu = \rho_2$, then $\gamma:= (X_1, X_2) \# d\mu$ has marginals $\rho_1$ and $\rho_2$, which implies that
$$
\int_{\Omega_0}d(X_1, X_2)^2 d\mu = \int_{\Omega \times \bR^3} d(x,y)^2 d\gamma (x,y) \leq W_2^2 (\rho_1, \rho_2).
$$
\end{Rmk}

\section{Proof of Theorem~\ref{theo1}}
\label{sec-prooftheo1}
In order to prove Theorem~\ref{theo1}, we shall focus on the case $\Omega =\T^3$ as  the treatment of space mean values requires a specific analysis compared to the case of $\R^3$.
We  explain in a final subsection the (slight) required modifications to handle the case $\Omega = \R^3$.

\subsection{Lagrangian formulation for weak solutions of Vlasov-Poisson and relativistic Vlasov-Maxwell}

We adopt a Lagrangian point of view, which means that our analysis will essentially rely on following the particles along their path. It means in concrete terms that we consider the two characteristic systems of ODEs corresponding to the Vlasov-Poisson system and the Vlasov-Maxwell system, starting at $(x,\xi)$ at time $0$.\\

We consider a weak solution $\fvp$ to \eqref{vp} in the sense of Definition \ref{defsolweakvp} and a weak solution $\fvm$ to \eqref{vm} in the sense of Definition \ref{defsolweakvm}.\\

The macroscopic densities $\rho^i$ are assumed to be bounded in $L^1 \cap L^{\infty}$ and therefore $\nabla_x \phi^i$ classically satisfies the following regularity properties (see  Lemma $3.1$ of \cite{Loeper-2006} or Lemma $3.2$ of \cite{Han-KwanIcaobelli-JDE}):
\begin{Lem}
\label{lem-loglip}
Let $\phi$ satisfy the Poisson equation
$$
- \Delta \phi = \rho-1, \quad \text{in }\, \bT^3.
$$
Then there exists $C$ depending only on $\|\rho-1\|_{L^{\infty}(0,T \times \bT^3)}$, such that
$$
\|\nabla \phi\|_{L^{\infty}(0,T \times \bT^3)} \leq C, 
$$
and
\begin{multline*}
\forall t \in [0,T), \, \forall (x,y)\in \bR^3 \times \bR^3, \\
|\nabla\phi(t,x) - \nabla\phi(t,y)| \leq C |x-y|\left(1+\log^+\left(\frac{1}{|x-y|}\right)\right),
\end{multline*}
where $\log^+ (z) = \log z$ if $z \geq 1$, $\log^+ (z) =0$ if $z < 1$.
\end{Lem}
This is enough to define a unique H\"older continuous flow (see e.g. \cite{Majda-Bertozzi})
that satisfies
 \be\label{charvp}
\left\{
\ba
{}&\d_t \Xvp (t,x,\xi) = \XiVp(t,x,\xi) , \\
& \d_t \XiVp (t,x,\xi) =-\nabla \phivp (t, \Xvp),  
\ea
\right.
\ee
with the initial conditions
\be
\left\{
\ba
{}& \Xvp(0,x,\xi)=x, \\
& \XiVp(0,x,\xi) = \xi.
\ea
\right.
\ee

On the other hand, for the Vlasov-Maxwell system, by definition of a weak solution, we also have  that the electromagnetic field $(E,B)$ belong to $L^{1}(0,T; BV(\bT^3)) \cap L^1(0,T; L^\infty(\bT^3))$. It follows from Ambrosio  \cite[Theorem 6.2]{Ambrosio-2004-Inventionnes}, that the solutions of the following characteristic systems of ODE exist, belong to $L^1_{loc}([0,T)\times\bT^3\times\bR^3)$ and  are absolutely continuous  (in time) for a.e $(x,\xi) \in \bT^3 \times \bR^3$:
\be\label{charvm} 
\left\{
\ba
{}& \d_t \Xvm(t,x,\xi) = v(\XiVm(t,x,\xi)) , \\
& \d_t \XiVm (t,x,\xi) =E(t, \Xvm) + \eps v(\XiVm)\times B(t,\Xvm),
\ea
\right.
\ee
endowed with the initial data
\be
\left\{
\ba
{}& \Xvm(0,x,\xi)=x, \\
& \XiVm(0,x,\xi) = \xi.
\ea
\right.
\ee
The solutions $(\Xvm, \XiVm)$ are unique accorded to \cite[Theorem 6.4]{Ambrosio-2004-Inventionnes} and form the Lagrangian flow associated to the Lorentz force field. We refer also to \cite{Ambrosio-Crippa-2008} for a review of further developments in this theory.\\

Moreover, we have the following representation formula:
\begin{equation}
  \forall t \in [0,T), \quad  \fvp= (\Xvp,\XiVp)(t,\cdot,\cdot)\# f^0,
\end{equation}
and
\begin{equation}
  \forall t \in [0,T), \quad    \fvm = (\Xvm,\XiVm)(t,\cdot,\cdot)\#f^0.
\end{equation}
Likewise,
\begin{equation}\label{rovp}
 \forall t \in [0,T), \quad  \rhoVp= \Xvp(t,\cdot,\cdot)\# f^0,
\end{equation}
and
\begin{equation}\label{rhovm}
 \forall t \in [0,T), \quad  \rhovm= \Xvm(t,\cdot,\cdot)\# f^0.
\end{equation}

\subsection {Log-Gr\"onwall estimate on the square of the $W_2$ distance between the Lagrangian trajectories}
%We consider an initial distribution $f^0(x,\xi)$ such that the corresponding weak solutions $\fvp$ and $\fvm$ to respectively \eqref{vp} and \eqref{vm} satisfy the conditions of the Theorem \ref{theo1}. \\
We  define the functional
\begin{equation}\label{Q}
Q(t) : = \frac{1}{2}\int_{\bT^3\times\bR^3}f^0(dx, d\xi)\left(|\Xvp-\Xvm|^2 + |\XiVp-\XiVm|^2\right)
\end{equation}
which quantifies the distance between the two solutions in a weak sense that we are going to explain.
One can notice that 
$$
((\Xvp,\XiVp)(t), (\Xvm,\XiVm)(t))\# f^0 \
$$
is a probability measure on $\left(\bT^3\times\bR^3\right)^2$ with marginals $\fvp$ and $\fvm$, which leads to the important preliminary lemma:
\begin{Lem}\label{W2Q}
Let $Q$ be the quantity defined in \eqref{Q}, then
$$
W_2^2(\fvp(t),\fvm(t)) \leq 2Q(t),
$$
and
$$
W_2^2(\rhoVp(t),\rhovm(t)) \leq 2Q(t).
$$
\end{Lem}

Any control of the functional $Q(t)$ will consequently imply an estimate of the Wasserstein distance between $\fvp$ and $\fvm$.

One can notice that the same considerations as above on the Lagrangian flows for the quantities $|\Xvp-\Xvm|^2$ and $|\XiVp-\XiVm|^2$ lead to
$$
|\Xvp-\Xvm|^2 (t)= 2 \int_0^t (\Xvp-\Xvm)(s)\cdot\left(\XiVp-v(\XiVm)\right)(s) ds,
$$ 
and
$$
|\XiVp-\XiVm|^2 (t)= 2 \int_0^t (\XiVp-\XiVm)(s)\cdot(F^{\text{VP}}(s,\Xvp)-F^{\text{VM}}(s,\Xvm)) ds.
$$
with the notation
$$
\ba
F^{\text{VP}} (t,x) &= - \nabla_x \phivp(t,x), \\
F^{\text{VM}} (t,x,\xi) &=- \nabla_x \phivm(t,x)- \eps \d_t A (t,x)+  \eps v(\xi)\times B(t,x).
\ea
$$
This yields for all $t$ in $[0,T)$:
\begin{equation}
\begin{aligned}
    \label{Qderivation}
 Q(t)&\leq  Q(0)\\
 &+\int_0^t\int_{\bT^3\times\bR^3}f^0(dx, d\xi) \left| \left(\Xvp(s)-\Xvm(s)\right)\cdot\left(\XiVp(s)-v(\XiVm(s))\right) \right|\,   ds \\
    &+ \int_0^t\int_{\bT^3\times\bR^3}f^0(dx, d\xi)  \Big| \left(\XiVp(s,x,\xi)-\XiVm(s,x,\xi)\right) \\
    &\qquad \qquad \qquad \cdot (\nabla_x\phivp(s,\Xvp(s,x,\xi))-\nabla_x\phivm(s,\Xvm(s,x,\xi))) \Big| \,   ds \\
    &+ \eps \int_0^t \int_{\bT^3\times\bR^3}f^0(dx, d\xi) \left| \left(\XiVp-\XiVm\right)\cdot(v(\XiVm)\times B(s,\Xvm)) \right| \,   ds \\
    &+ \eps \int_{\bT^3\times\bR^3}f^0(dx, d\xi) \left|\int_0^t\left(\XiVp-\XiVm\right)\cdot\d_t A(s,\Xvm)  \, ds \right|  .
\end{aligned}
\end{equation}
%\begin{Rmk}
%Under the hypothesis of the Theorem \ref{theo1}, $Q(0)=0$. But if we take the initial data depending on the parameter $\eps$ as explained in the Remark \ref{modificationspossibles}, then $Q(0)$ becomes proportional to $W_2^2(f^0, f^0_{\eps})$. It would change the rate of convergence by a simple modification of the following computations.
%\end{Rmk}

In the following, the $C$ will stand for a generic positive constant (independent of $\eps$ but that may depend on the initial data) that may change from line to line. \\

We will be able to somehow systematically replace $\XiVp(t)-v(\XiVm)(t)$ by $\XiVp(t) - \XiVm(t)$, up to some error terms. This is the content of the next lemma.

\begin{Lem}
\label{lem-vrelat}
 Let $G \in L^1_{loc}(\T^3 \times \R^3)$. For almost all $s \in [0,T)$, we have the estimate
\begin{multline}
\int_{\bT^3\times\bR^3}f^0(dx, d\xi)\left|\left(\XiVp-v(\XiVm)\right)\cdot G \right|       \\
\leq C\Bigg(  \eps^2 + \int_{\T^3\times \R^3}f^0(dx, d\xi) |G|^2 + \left |\int_{\bT^3\times\bR^3}f^0(dx, d\xi)\left|\XiVp-\XiVm\right| | G|  \right|   \Bigg).
\end{multline}

\end{Lem}

\begin{proof}
We can  write 
$$
\XiVp-v(\XiVm)= \left[\XiVp - v(\XiVp)\right] + \left[v(\XiVp)- v(\XiVm)\right]
$$
We observe then that 
$$
\left|v(\xi) - \xi \right| \leq \frac{\eps |\xi|^2}{\sqrt{1+\eps^2 |\xi|^2}} \leq \eps |\xi|^2,
$$
and by the Cauchy-Schwarz inequality and the Young inequality, we infer
\begin{multline*}
\int_{\bT^3\times\bR^3}f^0(dx, d\xi)\left| \left(\XiVp-v(\XiVp)\right)\cdot G(x,\xi)   \right|
\\
\leq \frac{1}{2}\left(\eps^2 \left |\int_{\T^3 \times \bR^3}|\xi|^4 \fvp(t,dx,d\xi)  \right| 
 +    \int_{\T^3\times \R^3}f^0(dx, d\xi) |G(x,\xi)|^2 \, d \xi dx\right),
\end{multline*}
and the first term is bounded by $C \eps^2$ thanks to the assumption that the solution to Vlasov-Poisson is suitable in the sense of Definition~\ref{sol-vp}, which implies
$$
\left\|\int_{\bR^3}|\xi|^4 \fvp(t,x, d\xi)\right\|_{L^{\infty}(0,T; L^1(\bT^3))} \leq C_0.
$$
Moreover, a straightforward computation ensures that the gradient of the velocity is bounded by a constant independent of $\eps$:
$$
\| \nabla v\|_{L^{\infty}(\bR^3)} \leq C.
$$
Therefore we have the estimate
$$
|v(\XiVp)- v(\XiVm)| \leq C |\XiVp - \XiVm|
$$
and we can conclude.
\end{proof}
We apply Lemma~\ref{lem-vrelat} to the first term in the expansion of the rhs of ~\eqref{Qderivation}, for $G = \Xvp(s)-\Xvm(s)$. We deduce a control by
$$ C\Bigg(  T \eps^2 + \int_0^t Q(s) \, ds    \Bigg).
$$
To derive a suitable estimate for $Q$ we therefore focus on the remaining terms of~\eqref{Qderivation}. We define the following three quantities $I_1, I_2$ and $I_3$  for $t \in [0,T)$, which we will  tackle one  after another:
\begin{align*}
I_1&:=  \int_{\bT^3\times\bR^3}f^0(dx, d\xi) \left|\Xvp(t)-\Xvm(t)\right| \left|\XiVp(t)-\XiVm(t)\right|\\
&+ \int_{\bT^3\times\bR^3}f^0(dx, d\xi)\Big| \left(\XiVp(t)-\XiVm(t)\right) \\
&\qquad \qquad \qquad \cdot(\nabla_x\phivp(t,\XiVp(t))-\nabla_x\phivm(t, \XiVm(t)))\Big|,
\end{align*}
which will be estimated following the path traced in Section $3.2$ of \cite{Loeper-2006},
$$
I_2 :=  \eps\int_{\bT^3\times\bR^3}f^0(dx, d\xi)\left| \left(\XiVp(t)-\XiVm(t)\right)\cdot(v(\XiVm)\times B(t,\Xvm))\right|,
$$
whose estimate is almost straightforward with the uniform control~\eqref{growth-fields} on the $L^{2}$ norm of the magnetic field $B$, and
\be
\label{def-I3}
I_3 := \eps \int_{\bT^3\times\bR^3} f^0(dx, d\xi) \left|\int_0^t \left(\XiVp(t)-\XiVm(t)\right)\cdot\d_t A(t,\Xvm(t)) \, ds \right|,
\ee
which requests a little more subtle integration by part arguments, which explains that we need to keep the integral in time for estimating this contribution.%consider both integral in the phase space and in time together.

\subsubsection{Estimate on $I_1$}
In this paragraph we follow carefully the steps of the proof of Section $3$ of \cite{Loeper-2006}.\\

By a straightforward Young inequality the first term of $I_1$ is bounded by $Q(t)$. The Cauchy-Schwarz inequality implies that the second term is bounded by:
\begin{multline*}
\left(2Q(t)\right)^{1/2} \times \\
\left(\int_{\bT^3 \times\bR^3} f^0(dx, d\xi) \left|\nabla_x \phivp(t,\Xvp(t,x,\xi)) - \nabla_x \phivm(t,\Xvm(t,x,\xi))\right|^2 \right)^{1/2}.
\end{multline*}
We then write
\begin{align*}
& \left(\int_{\bT^3 \times\bR^3} f^0(dx, d\xi) \left|\nabla_x \phivp(t,\Xvp(t,x,\xi)) - \nabla_x \phivm(t,\Xvm(t,x,\xi))\right|^2 \right)^{1/2} \\
&\leq  \left(\int_{\bT^3 \times\bR^3} f^0(dx, d\xi) \left|\nabla_x \phivp(t,\Xvm(t,x,\xi)) - \nabla_x \phivm(t,\Xvm(t,x,\xi))\right|^2 \right)^{1/2} \\
&+  \left(\int_{\bT^3 \times\bR^3} f^0(dx, d\xi) \left|\nabla_x \phivp(t,\Xvp(t,x,\xi)) - \nabla_x \phivp(t,\Xvm(t,x,\xi))\right|^2 \right)^{1/2} \\
&=:J_1(t)^{1/2} + J_2(t) ^{1/2}.
\end{align*}
We are now going to estimate $J_1$ and $J_2$.\\

For $J_1$, the equation \eqref{rovp} and Theorem \ref{theo3}  imply that
\begin{align*}
J_1(t) &= \int_{\bT^3} \rhoVp(t,x) |\nabla_x \phivp(t,x) - \nabla_x \phivm(t,x)|^2 \, dx  \\ 
 &\leq \max \{\|\rhoVp\|_{L^{\infty}}, \|\rhovm\|_{L^{\infty}}\}^2  W_2^2 (\rhoVp(t), \rhovm(t)).
\end{align*}
Using  Lemma \ref{W2Q}, we conclude that
$$
J_1(t) \leq 2 \max \{\|\rhoVp\|_{L^{\infty}}, \|\rhovm\|_{L^{\infty}}\}^2 Q(t).
$$
By assumption both $\|\rhoVp\|_{L^{\infty}}$ and $\|\rhovm\|_{L^{\infty}}$ are uniformly bounded in $\eps$ (see \eqref{bound-rho-vm}).
The estimate of $J_2$ can be done from standard arguments relying on the Log-Lipschitz regularity of $\nabla_x \phivp$, see Lemma~\ref{lem-loglip}. We refer to the end of the section $3$ of \cite{Loeper-2006} for the computations leading to the following estimate:
$$
J_2(t) \leq C Q(t)\left(1+ \log^+\left(\frac{1}{Q(t)}\right)\right).
$$
Gathering the previous estimates finally gives:
\begin{equation}
\label{I1}
I_1 \leq C Q(t) \left(1+ \log^+\left(\frac{1}{Q(t)}\right)\right).
\end{equation}

\subsubsection{Estimate on $I_2$}
First, an application of the Cauchy-Schwarz inequality gives:
\begin{align*}
I_2 &\leq  \eps Q(t)^{1/2}\left(\int_{\bT^3\times\bR^3}f^0(dx,d\xi)|v(\XiVm)|^2 |B(t,\Xvm)|^2 \right)^{1/2}\\
          &\leq \eps Q(t)^{1/2}\left(\int_{\bT^3\times\bR^3}\fvm(t,dx,d\xi)|v(\xi)|^2|B(t,x)|^2 \right)^{1/2} \\
          &\leq \eps Q(t)^{1/2} \|B(t)\|_{L^2(\bT^3)} \left\|\int_{\bR^3}  \fvm(t,x,d\xi) |v(\xi)|^2  \right\|^{1/2}_{L^{\infty}(\bT^3)}.
\end{align*}
By the assumption~\eqref{growth-fields}, we have 
$$
\|B(t)\|_{L^2(\bT^3)} \leq C\eps^{-\gamma_2}.
$$
It remains then to estimate the quantity $\eps^2\|\int_{\bR^3} \fvm(t,x,d\xi) |v(\xi)|^2\|_{L^{\infty}(\bT^3)}$.\\

To this end, we use the uniform bound~\eqref{malpha} bearing on
 $\|m_\alpha\|_{L^{\infty}(\bT^3)}$ and  the general fact that $\eps|v(\xi)|\leq 1$. This gives 
$$
\ba
\eps^2\left\|\int_{\bR^3}  \fvm(t,x,d\xi) |v(\xi)|^2 \right\|_{L^{\infty}(\bT^3)}&\leq C \eps^\alpha \left\|\int_{\bR^3} \fvm(t,x,d\xi)  |v(\xi)|^{\alpha} \right\|_{L^{\infty}(\bT^3)} \\
&\leq C\eps^{\alpha-\beta}.
\ea
$$
Consequently we obtain using Young's inequality,
\begin{equation}
\label{I2}
I_2 \leq C \eps^{\alpha-\left(\beta+2\gamma_2\right)} + C Q(t).
\end{equation}

\subsubsection{Estimate on $I_3$}
\label{sec-I3}
The same direct approach fails for the estimate of $I_3$ because it gives
$$
I_3 \leq C Q(t)^{1/2} \|\eps \d_t A(t,\Xvp(t,.))\|_{L^2(\bT^3)},
$$
but unfortunately, at first glance, we can only use~\eqref{growth-fields}, that gives
$$
\|\eps \d_t A\|_{L^2(\bT^3)} \leq C\eps^{-\gamma_1},
$$
and we seemingly lose some power of $\eps$.

\begin{Rmk}
Even if we assume some uniform bound in $L^2$ for $\jvm$,
the same problem of "loss of $\eps$" will occur. Indeed, the energy estimate for the wave equation
$$
\eps^2\partial_t^2A-\Delta_x A=\eps \bP(\jvm)
$$
gives the estimate
$$
\|\eps \d_t A(t,\Xvp(t,.))\|_{L^2(\bT^3)}\leq CT\|j\|_{L^\infty(0,T;L^2(\bT^3))},
$$
and the small parameter $\eps$ is lost as well.\\
\end{Rmk}
This is the reason why we have to deal with the Gronwall inequality in its integral form in order to perform an integration by parts with respect to the time variable.\\
One must first observe that
\begin{align*}
\d_s A(s,\Xvm(s,x,\xi))=& \frac{d}{ds}\left(A\left(s,\Xvm(s,x,\xi)\right)\right) - \d_s \Xvm \cdot \nabla_x A(s,\Xvm)\\
            =& \frac{d}{ds}\left(A\left(s,\Xvm(s,x,\xi)\right)\right) - v( \XiVm)\cdot \nabla_x A(s,\Xvm).
\end{align*}
We then define
$$
I_{31} := \eps\int_{\bT^3\times\bR^3}f^0(dx, d\xi) \left| \ \int_0^t   \frac{d}{ds}\left[A\left(s,\Xvm(s,x,\xi)\right)-\langle A \rangle \right] \cdot (\XiVp-\XiVm) ds \right|,
$$
$$
I_{32} :=  \eps\int_{\bT^3\times\bR^3} \int_0^t f^0(dx, d\xi) \left| \left(v(\XiVm)\cdot \nabla_x A(s,\Xvm)\right) \cdot (\XiVp-\XiVm)\right| ds,
$$
and
$$
I_{33} :=  \eps\int_{\bT^3\times\bR^3} \int_0^t f^0(dx, d\xi) \left| \left(\langle \d_s A \rangle \right)\cdot  (\XiVp-\XiVm)(s) \right| ds,
$$
so that
$$
I_3 \leq I_{31} + I_{32} + I_{33}.
$$
To estimate $I_{32}$ we perform again a Cauchy-Schwarz inequality in the integral over $\bT^3$:
$$
I_{32} \leq \eps \int_0^t \left(Q(s)^{1/2}\left( \int_{\bT^3\times \R^3}|v(\xi)|^2\fvm(s,dx,d\xi) |\nabla_x A (s,x)|^2 \right)^{1/2} \right)ds.
$$
The factor $\nabla_x A$ can be bounded in $L^2$ by the $L^2$ norm of $B$  thanks to the Biot and Savart law.
\begin{Lem}\label{biotetsavart}
Let $B\in L^2(\bT^3)$ and $A$ such that
$$
\nabla\times A = B - \langle B \rangle.
$$
Then we have the Biot and Savart law:
$$
A -\langle A \rangle = \nabla\times \left(\Delta^{-1}\left(B-\langle B\rangle\right)\right).
$$
where $\Delta^{-1}$ selects the unique solution with zero mean to the associated Poisson equation.
It follows in particular that
$$
\|\nabla A\|_{L^2(\bT^3)}\leq C \|B - \langle B \rangle\|_{L^2(\bT^3)}.
$$
\end{Lem}
\begin{proof}
There exist a unique distribution $\psi$ such that
$$
-\Delta \psi = B - \langle B \rangle, \quad \int_{\T^3} \psi \, dx =0.
$$
Then one can check that
$$
\nabla \times \left( A - \nabla\times \psi \right) = - \nabla\left(\nabla\cdot \psi \right) =0,
$$
and therefore
$$
A = \nabla \times \psi  + \langle A \rangle.
$$
The Biot and Savart law and the estimate follow.
\end{proof}

We also have the following conservation of the spatial mean of $B$.

\begin{Lem}
The space mean-value of $B$ is constant, for any $t$ in $[0,T)$:
$$
\langle B(t) \rangle = \langle B^0 \rangle
$$
\end{Lem}

\begin{proof}
It is straightforward since from the Maxwell equations, for any $t$ in $[0,T)$ we have
$$
\eps \frac{d}{dt}\langle B(t) \rangle = \eps \langle \d_t B(t) \rangle = - \langle \nabla\times E(t) \rangle = 0.
$$
\end{proof}

We deduce from the above lemmas and~\eqref{growth-fields} that
$$
\|B - \langle B \rangle\|_{L^2(\bT^3)} \leq \|B\|_{L^2(\bT^3)} + |\langle B \rangle| \leq C\eps^{-\gamma_2},
$$
so that
$$
\|\nabla_x A\|_{L^{\infty}(0,T;L^2(\bT^3))}\leq C\eps^{-\gamma_2}.
$$
We  then have using the Young inequality,
\begin{equation}
\label{I32}
I_{32} \leq C \int_0^t Q(s) ds + C \eps^{\alpha -(\beta+\gamma_2)}.
\end{equation}
To estimate $I_{33}$ we first use the Cauchy-Schwarz inequality
$$
I_{33} \leq \eps\int_0^t Q(s)^{1/2}|\langle \d_s A(s)\rangle| ds
$$
and then  rely on the fact that for any $t$ in $[0,T)$, since $\eps \langle \d_t A \rangle$ satisfies
$$
\eps \frac{d}{dt} \langle \d_t A \rangle =  \langle \jvm \rangle,
$$
we have
\begin{equation}
\label{moy-A}
\eps \langle \d_t A \rangle = \int_0^t \langle \jvm \rangle (s) ds + \eps \langle \d_t A \rangle|_{t=0}.
\end{equation}
The initial data being normalized, the last term is by assumption on the initial electric field $E_\eps^0$ equal to $0$.
We therefore must focus on the first term of~\eqref{moy-A}.\\

We write that
$$
 \int_0^t \langle \jvm \rangle (s) ds =  \int_0^t \left(\langle \jvm \rangle (s) - \langle \jvp \rangle (s)\right)ds + \int_0^t \langle \jvp \rangle (s) ds.
$$
We then remark that
\begin{align*}
\int_0^t \left(\langle \jvm \rangle (s) - \langle \jvp \rangle (s)\right)ds &= \int_0^t \int_{\bT^3\times \bR^3}f^0(dx, d\xi) \left(v(\XiVm) - \XiVp\right) \,ds\\
&= \int_0^t \int_{\bT^3\times\bR^3}f^0(dx, d\xi) \left(v(\XiVm)- v(\XiVp)\right)\,  ds \\
&+ \int_0^t \int_{\bT^3\times\bR^3}f^0(dx, d\xi) \left(v(\XiVp)- \XiVp\right)\,  ds ,
\end{align*}
and since
$$
|v(\xi)- \xi| \leq \eps |\xi|^2
$$
and
$$
\|\nabla_{\xi}v\|_{L^{\infty}} \leq C,
$$
we have
$$
\left|\int_0^t \left(\langle \jvm \rangle (s) - \langle \jvp \rangle (s)\right)ds\right| \leq CT Q(t)^{1/2} + C\eps T\left\|\int_{\bR^3}|\xi|^2 \fvp (t,x,d \xi)\right\|_{L^1(\bT^3)}.
$$
Concerning $\int_0^t \langle \jvp \rangle (s) ds$, we use the fact that the Vlasov-Poisson equation preserves the current density which therefore will be equal to $0$ because  the initial data is normalized so that $\langle \jvp\rangle|_{t=0}=0$.

\begin{Lem}
We have the conservation of the spatial mean-value of the current density for the Vlasov-Poisson system:
$$
\forall t \in [0,T),\, \langle \jvp \rangle (t) =  \langle \jvp \rangle (0) .
$$
\end{Lem}
\begin{proof}
We have the following conservation law ensured by Definition \ref{defsolweakvp}
$$
\d_t \jvp + \nabla_x : \langle f ,1\otimes |\xi|^2 \rangle = \rhoVp E.
$$
Therefore
$$
\frac{d}{dt}\langle \jvp \rangle = \int_{\bT^3}\Delta_x \phivp \nabla_x \phivp \, dx.
$$
For any $i\in\{1,2,3\}$, we have
$$
\int_{\bT^3}\d_{ii} \phivp \d_i\phivp \, dx = \frac{1}{2}\int_{\bT^3}\d_i|\d_i\phivp|^2 \, dx = 0.
$$
Likewise, for any $i\neq j$, we have by integration by parts with respect to $x_i$,
\begin{align*}
\int_{\bT^3}\d_{ii}\phivp \d_j\phivp \, dx & =- \int_{\bT^3}\d_i\phivp \d_{ij}\phivp \, dx  \\
&=-\frac{1}{2}\int_{\bT^3}\d_j|\d_i\phivp|^2 \, dx \\
&= 0.
\end{align*}
The lemma is finally proved.
\end{proof}

We end up with the following estimate for $I_{33}$, for any $t$ in $[0,T)$:
\begin{equation}
\label{I33}
 I_{33}(t) \leq C (1+T)^2 \left(\eps^2+ \int_0^t  Q(s)ds\right) .   
\end{equation}

To estimate $I_{31}$ we first perform the integration by parts with respect to the time variable, which yields
\begin{equation}
\label{I31-bound}
|I_{31}|\leq  K_1 + K_2,
\end{equation}
with
\begin{align*}
&K_1 := \eps\int_{\bT^3\times\bR^3} \int_0^t f^0(dx, d\xi)\\
&\qquad \qquad\qquad\left| (A(s,\Xvm(s,x,\xi))- \langle A \rangle(s))\cdot\left(F^{\text{VP}}(\Xvp(t))-F^{\text{VM}}(\Xvm(t))\right) \right|ds, \\
&K_2 := \eps\int_{\bT^3\times\bR^3} f^0(dx, d\xi) \left| \left(A(t,\Xvm(t))- \langle A \rangle(t)\right)\cdot\left(\XiVp(t) - \XiVm(t)\right) \right|,
\end{align*}
where we recall
\begin{multline*}
F^{\text{VM}}(s,\Xvm)-F^{\text{VP}}(s,\Xvp) =\\
 \nabla_x \phivm(s,\Xvm) - \nabla_x \phivp(s,\Xvp) + \eps v(\XiVm)\times B(s,\Xvm) + \eps \d_t A(s,\Xvm).
\end{multline*}
We first treat $K_1$. We are somehow back to the terms $I_1$, $I_2$  but with a gain of a power of $\eps$. Then by~\eqref{growth-fields},
$$
\|\eps \d_t A \|_{L^2(\bT^3)} (t) \leq C_0\eps^{-\gamma_1},
$$
and performing the same analysis, one has
\begin{align*}
|K_1|\leq &  C\|\rhovm\|_{L^{\infty}(\bT^3)}\|A(t,.)- \langle A \rangle\|_{L^{2}(\bT^3)}\\
& \times \left(\eps \int_0^t  Q(s)\left(1+ \log^+\left(\frac{1}Q(s)\right)\right) ds + \eps^{\alpha-(\beta+\gamma_2)}   +  \eps^{1-\gamma_1}T  \right).
\end{align*}
Now by the Poincar\'e inequality on $\T^3$ and the Biot and Savart law, we have
\begin{equation}
\label{Poincare}
\|A(t,.)- \langle A \rangle\|_{L^{2}(\bT^3)} \leq C \| \nabla A \|_{L^{2}(\bT^3)}  \leq C \| B- \langle B\rangle\|_{L^2(\T^3)}  \leq C \eps^{-\gamma_2}.
\end{equation}

Summing up all these estimates,  
we obtain that for all $t$ in $[0,T)$:
\begin{equation}
\label{K1}
|K_1|\leq C(1+T)^2\left( \eps^{\alpha-(\beta + 2 \gamma_2)} + \eps^{1-(\gamma_1+\gamma_2)} +  \int_0^t Q(s)\left(1+ \log^+\left(\frac{1}Q(s)\right) \right)ds \right).
\end{equation}
For $K_2$, we get by Cauchy Schwarz
$$
|K_2| \leq \eps \| \rhovm(t)\|_{L^\infty} \|( A- \langle A \rangle)(t)\|_{L^2(\bT^3)}  Q(t)^{1/2}.
$$
Therefore, by the Young inequality, we conclude that
\begin{equation}
\label{K2}
\ba
|K_2| &\leq \eps^{1+\gamma_2 - \gamma_1} \| \rhovm(t)\|_{L^\infty}^2  \|( A- \langle A \rangle)(t)\|_{L^2(\bT^3)}^2 + \eps^{1+\gamma_1- \gamma_2} Q(t) \\
&\leq C \eps^{1- (\gamma_1+ \gamma_2)} +  C \eps^{1+\gamma_1- \gamma_2} Q(t),
\ea
\end{equation}
where we have used again the Biot and Savart law and~\eqref{growth-fields}.
Finally, gathering~\eqref{I32}, \eqref{I33}, \eqref{I31-bound}, \eqref{K1} and \eqref{K2}, since $\gamma_2<1$, we obtain
\begin{multline}
\label{I3}
|I_3| \leq C(1+T)^2 \eps^{\min(\alpha- (\beta+2 \gamma_2), 1-(\gamma_1+\gamma_2))} \\
+ \int_0^t C(1+T)^2 Q(s) \left(1+ \log^+\left(\frac{1}{Q(s)}\right)\right) ds +  C \eps^{1+\gamma_1- \gamma_2}  Q(t).
\end{multline}

\subsection{Final estimate}
Finally by \eqref{I1}, \eqref{I2} and \eqref{I3} (for $\eps \in [0,\eps_0]$ with $\eps_0$ enough to absorb the term $ C \eps^{1+\gamma_1- \gamma_2}  Q(t)$ in~\eqref{I3} by $Q(t)$ of the left-hand side) we have the following Osgood estimate on the quantity $Q(t)$:
\begin{equation}
\label{eq-Osgood}
Q(t)\leq C\left(1+T\right)^2 \eps^{\kappa} + \int_0^t C\left(1+T\right)^2 Q(s)\left(1+ \log^+\left(\frac{1}{Q(s)}\right)  \right) ds.
\end{equation}
with 
$$
\kappa= \min(\alpha- (\beta+2 \gamma_2),1-(\gamma_1+\gamma_2)).
$$
The procedure to obtain~\eqref{inegalitenonrelativiste} from~\eqref{eq-Osgood} is standard: let us quickly explain it for the sake of completeness. Set
$$
\ba
\mu(z) &= C\left(1+T\right)^2 z (1+ \log^+ (1/z)), \\
\varphi(t)&=C\left(1+T\right)^2 \eps^{\kappa} + \int_0^t C\left(1+T\right)^2 Q(s)\left(1+ \log^+\left(\frac{1}{Q(s)}\right)  \right) ds.
\ea
$$
Since $\mu$ is non-decreasing, we have
$$
\varphi'(t) = \mu(Q(t)) \leq \mu(\varphi(t)).
$$
Set then
$$
U(t)=  \log^+ \varphi(t).
$$
It follows that $U$ satisfies
$$
U'(t) \leq C(1+T)^2 (1- U(t)).
$$
which we can explicitly integrate, yielding
$$
U(t) \leq U(0)e^{-C(1+T)^2 t} + (1-e^{-C(1+T)^2 t}).
$$
Coming back to $Q$, by a continuity (in time) argument, taking $\eps_0>0$ small enough, we finally obtain that for all $\eps \in(0, \eps_0]$ and $t \in [0,T]$,
$$
Q(t)\leq C \exp\left(\log\left(C(1+T)^2\eps^{\kappa}\right)\exp\left(-C(1+T)^2\right)\right),
$$
which implies the desired inequality \eqref{inegalitenonrelativiste} and the proof of Theorem~\ref{theo1} is complete.

\subsection{The case $\Omega= \R^3$}
As already mentioned, the proof for $\Omega=\bR^3$ is very similar, yet simplified in some aspects. The main difference is that we do not need to handle space mean values as in the torus case.
This is in particular apparent in the treatment of the term $I_3$ (as defined in~\eqref{def-I3}).
We have this time
$$
I_3 \leq I_{31} + I_{32},
$$
where
$$
I_{31} = \eps\int_{\bR^3\times\bR^3} f^0(dx, d\xi) \left| \int_0^t  \frac{d}{ds}\left[A\left(s,\Xvm(s,x,\xi)\right) \right] \cdot (\XiVp-\XiVm) \, ds \right|,
$$
$$
I_{32} =  \eps\int_{\bR^3\times\bR^3} \int_0^t f^0(dx, d\xi) \left| \left(v(\XiVm)\cdot \nabla_x A(s,\Xvm)\right) \cdot (\XiVp-\XiVm) \right| ds  ,
$$
To study $I_{31}$, we rely on the same integration by parts in time argument. Only the final estimate is different: in $\bR^3$ the Biot and Savart law gets simplified compared to the case of $\bT^3$,
so that
$$
\| \nabla_x A \|_{L^2(\R^3)}\leq C \| B \|_{L^2(\R^3)}
$$
and we use the Sobolev embedding instead of the Poincar\'e inequality \eqref{Poincare}, which yields
$$
\ba
\| A (s) \|_{L^6(\R^3)} &\leq C \| \nabla_x A \|_{L^2(\R^3)} \\
&\leq C \| B (s) \|_{L^2(\R^3)}.
\ea
$$
Writing
$$
|I_{31}|\leq  K_1 + K_2,
$$
with
\begin{align*}
&K_1 := \eps\int_{\bR^3\times\bR^3} \int_0^t f^0(dx, d\xi)\\
&\qquad \qquad\qquad\left| (A(s,\Xvm(s,x,\xi))\cdot\left(F^{\text{VP}}(\Xvp)-F^{\text{VM}}(\Xvm)\right) \right|ds , \\
&K_2 := \eps\int_{\bR^3\times\bR^3} f^0(dx, d\xi) \left| (A(t,\Xvm(t))\cdot\left(\XiVp(t) - \XiVm(t)\right) \right| ,
\end{align*}
 the outcome is the estimate
\begin{align*}
|K_1|\leq   C&\|\rhovm\|_{L^{3}(\bR^3)}\|B(t,.)\|_{L^{2}(\bR^3)}\\
& \times \left(\eps \int_0^t  Q(s)\left(1+ \log^+\left(\frac{1}Q(s\right)\right) ds + \eps^{\alpha-(\beta+\gamma_2)}   +  \eps^{1-\gamma_1}T  \right), \\
|K_2| \leq \eps &\| \rhovm(t)\|_{L^{6/5}(\bR^3)}  \| B (t)\|_{L^2(\bR^3)} + \eps Q(t).
\end{align*}
The remaining of the proof applies, \emph{mutatis mutandis}.

\section{A  class of measure-valued solutions which satisfies the assumptions of Theorem~\ref{theo1}}
\label{sec-multi}

Let $\Omega =\bT^3$.
The goal of this section is to build measure-valued solutions to~\eqref{vm} and \eqref{vp} that are not in the class of compactly supported $C^1$ solutions in $x$ and $\xi$, and to which Theorem~\ref{theo1} can nevertheless apply.  %for which the classical limit results have been derived in the literature \cite{Degond-86-MMAS}, \cite{Asano-Ukai-86-SMA}, \cite{Schaeffer-86-CMP}. 
More precisely we are interested in solutions with high regularity in the space variable $x$, namely real-analytic, and very little in the momentum variable $\xi$: basically we only ask for $f(t,x,.)$ to be a measure with some finite moments. The corresponding weak solutions to \eqref{vm} and \eqref{vp} will in fact be induced by a family of strong solutions in $x$ to a related fluid system.
In the first subsection we will give the definition of  what we call weak in $\xi$ and strong in $x$ solutions, following Baradat \cite{Aymeric-2019}. Then we explain a multifluid representation (used by Grenier \cite{Grenier-96}) that will allow us to rewrite the Vlasov-Maxwell system as a system of  fluid equations that we will effectively study.  Finally we will prove a small time existence result for these systems by a Cauchy-Kovalevskaya argument, again following \cite{Grenier-96} (see also \cite{GN}).

\subsection{Weak in $\xi$ and strong in $x$  solutions}
We will consider in the following a particular class of the more general solutions we handled previously in  Definitions \ref{defsolweakvp} and \ref{defsolweakvm}. It concerns 
weak solutions that are regular in $x$, at least $C^1$, for which another convenient definition can be given,  following Baradat \cite{Aymeric-2019}.

Let $T>0$.
We consider a function
$
f: [0,T] \times\bT^3 \rightarrow \mathcal{P}(\bR^3),
$
such that for any test function $\varphi \in C^1_b (\bR^3)$, the hydrodynamic observable corresponding to $\varphi$:
$$
\langle f,\varphi\rangle (t,x) := \int_{\bR^3}\varphi(\xi)f(t,x,\, d\xi),
$$
is a smooth function, namely in $C^1([0,T)\times \bT^3)$. 

\begin{Def}

We say that $f: [0,T] \times\bT^3 \rightarrow \mathcal{P}(\bR^3)$ is a weak in $\xi$ and strong in $x$  solution to \eqref{vp} if it satisfies in the classical sense, for all test functions $\varphi\in C^1_b(\bR^3)$ the system
\be
\left\{
\ba
{}&\d_t \langle f,\varphi\rangle (t,x) + \nabla_x\cdot \langle f,\xi\varphi\rangle (t,x) + \nabla_x \phivp \cdot \langle f, \nabla\varphi\rangle (t,x) =0, \\
&-\Delta_x \phivp = \langle f, 1\rangle  -1,\\
&f (0,x, \, d\xi) = f^0(x, \, d\xi).
\ea
\right.
\ee

\end{Def}

\begin{Def}

We say that $f: [0,T] \times\bT^3 \rightarrow \mathcal{P}(\bR^3)$ is a weak in $\xi$ and strong in $x$ solution to \eqref{vm} if it satisfies in the classical sense for all test functions $\varphi\in C^1_b(\bR^3)$ the system
\be \label{stronginxweakinv}
\left\{
\ba
{}&\d_t \langle f,\varphi\rangle (t,x) + \nabla_x\cdot \langle f,v(\xi)\varphi\rangle (t,x) \\
&\quad +E \cdot \langle f, \nabla \varphi\rangle + \eps \sum_{i=1}^3  B_{\sigma^2(i)} \langle  f,   v_{\sigma(i)}(\xi) \partial_{_i}\varphi \rangle   - B_{\sigma(i)} \langle f, v_{\sigma^2(i)}(\xi) \partial_{i}\varphi  \rangle =0, \\
&-\Delta_x \phi = \langle f, 1\rangle  -1,\\
&\eps^2 \d_t A - \Delta_x A = \eps \bP(\langle f, v(\xi)\rangle ), \\
&f(0,x, \, d\xi) = f^0(x, \, d\xi),\\
&\eps\d_t A|_{t=0} = E^0 + \nabla_x \phi|_{t=0},\\
&\nabla\times A|_{t=0} = B^0 - \langle B^0 \rangle, \quad \langle A|_{t=0} \rangle =0,
\ea
\right.
\ee
where $\sigma$ stands for the permutation 
$\sigma = (1,2,3)$.
\end{Def}

\subsection{A multifluid representation}
In this section we will set another formulation of the Vlasov-Maxwell system of equations, which we refer to as the multifluid representation, as introduced in \cite{Grenier-96} to prove a small time uniform existence result and analyze the quasineutral limit for Vlasov-Poisson type systems. It has also been used in \cite{Aymeric-2019} for studying nonlinear instabilities around rough velocity profiles in Vlasov-Poisson systems. We will be able to prove the existence of strong solutions to this system which in turn will provide some weak in $\xi$ and strong in $x$ solutions to the Vlasov-Maxwell system.\\

We look for solutions $f$ under the  form
\begin{equation}
\label{def-f-multi}
f(t,x,d \xi) = \int_{M}\rho^{\theta}(t,x)\delta(\xi - \xi^{\theta}(t,x)) \mu(d\theta),
\end{equation}
where $(M,\mu)$ is (fixed) a probability space, $\delta$ stands for the Dirac mass at $0$, $(\rhoteta)_{\theta\in M}$, $(\xiteta)_{\theta\in M}$ are families of smooth functions and vector fields on $\bT^3$. This is a representation where the whole set of particles in the plasma can be divided into several phases, each of them characterized by its pointwise macroscopic density $\rhoteta(t,x)$ and its pointwise momentum $\xiteta(t,x)$. Each density will be transported by the relativistic velocity $v(\xiteta)$ according to a continuity equation and each phase will be accelerated by the same electromagnetic field, producing a Lorentz force $F_L$ that is computed by taking into account all the different phases.\\

More precisely, given smooth initial data $(\rhoteta^0)_{\theta\in M}$, $(\xiteta^0)_{\theta\in M}$ we consider the following system:
\be\label{multiphasic}
\left\{
\ba
 &\forall \theta\in M,\, \d_t \rhoteta(t,x) +\nabla_x \cdot\left( v(\xiteta(t,x)) \rhoteta(t,x)\right)=0  \,,
\\
&\forall \theta\in M ,\, \d_t\xiteta (t,x) + \left(v(\xiteta(t,x))\cdot\nabla_x\right)\xiteta(t,x) \\&\qquad \qquad \qquad \qquad = \left( -\nabla_x \phi - \eps \d_t A +\eps v(\xiteta)\times (\nabla_x\times A + \langle B^0 \rangle\right)(t,x) \,,\\
&\eps^2\partial_t^2 A-\Delta_x A=\eps \bP \left(\int_{\bT^3}v(\xiteta(t,x))\rhoteta(t,x) \mu(d\theta)\right) \,,
\\  
&- \Delta_x \phi = \int_{M}\rhoteta(t,x) \mu(d\theta) - 1 \,,\\
&\forall \theta\in M, \quad \rhoteta|_{t=0}=\rhoteta(0), \quad \xiteta|_{t=0}=\xiteta(0), \\
&\eps\d_t A|_{t=0} = E^0 + \nabla_x \phi|_{t=0},\\
&\nabla\times A|_{t=0} = B^0 - \langle B^0 \rangle, \quad \langle A|_{t=0} \rangle =0,
\ea
\right.
\ee
Note that this corresponds to an initial condition 
$$
f^0(x,d\xi)= \int_{M}\rhoteta(0,x)\delta(\xi - \xiteta(0,x)) \mu(d\theta).
$$
As explained in \cite{Grenier-96}, this allows to model a variety of initial conditions, including
\begin{itemize}
\item continuous functions in $x$ and $\xi$, taking $(M, \mu)=\left(\R^3, \lambda \frac{d\theta}{1+\theta^2}\right)$ (where $\lambda>0$ is chosen in order to normalize the measure), 
$$\xiteta(0,x)= \theta, \quad \rhoteta(0,x) = \pi(1+\theta^2) f(t,x,\theta);$$
\item finite sums of Dirac masses in velocity supported on $v_1, \cdots, v_n \in \R^3$, that corresponds to a sum of monokinetic data, in which case $(M,\mu)$ is a discrete probability space with uniform measure.
\end{itemize}

Assuming that we are able to solve this system, we can define the measure $f(t,x,.)$ according to~\eqref{def-f-multi}. Given a smooth test function $\varphi$, we have then
$$
\langle f,\varphi\rangle(t,x) = \int_M \varphi(\xiteta(t,x))\rhoteta(t,x)\mu(d\theta).
$$
This is a straightforward computation to check that if $(\rhoteta)_{\theta\in M}$, $(\xiteta)_{\theta\in M}$ solve \eqref{multiphasic} then the measure valued function $f$ defined above is a strong in $x$ and weak in $\xi$ solution to \eqref{stronginxweakinv}.\\

As System~\eqref{multiphasic} does not seem to possess any hyperbolic structure, we are forced to solve it for initial data with analytic regularity, using a Cauchy-Kovalevskaya type scheme. The precise analytic spaces we work with are as follows.
\begin{Def}
For $\delta>1$, $B_{\delta}$ is the space of real functions $f$ on $\bT^3$ such that
$$
|f|_{\delta}:= \sum_{k\in\bZ^3}\left|\cF f (k)\right| \delta^{|k|}  < + \infty,
$$
where $(\cF f(k))$ are the Fourier coefficients of $f$. %(here we work on $\bT^3$, hence it is a function on $\bZ^3$).
\end{Def}  

We are now able to state the second main result of our paper.

\begin{Thm}
\label{thm2}
Let $(M,\mu)$ the probability space used to define the multifluid system \eqref{multiphasic}, let $\delta_0> \delta_1>1$, let $C_0>0$ and let $(\rhoteta(0))_{\theta}$, $(\xiteta(0))_{\theta}$ and $(E^0_\eps, B^0_\eps)$ be  families of $\bdeltazero$ satisfying
    $$
\nabla \cdot E^0_\eps = \int_{M} \rhoteta(0) \, d\mu(\theta), \qquad \nabla\cdot B^0_\eps = 0,
    $$
and such that
\begin{equation}
\begin{aligned}
&\underset{\theta}{\sup} |\rhoteta (0) |_{\delta_0} \leq C_0, \\
&\underset{\theta}{\sup} |\xiteta(0) |_{\delta_0} \leq C_0.
\end{aligned}
\end{equation}
Assume also that for some $\gamma \in [0,1]$,
\begin{equation}
\label{E0B0init}
|E^0_\eps|_{\delta_0} +  |B^0_\eps |_{\delta_0} \leq C_0 \eps^{-\gamma}.
\end{equation}
Then there exist a constant $\eps_0>0$, such that for any $\eps \in (0,\eps_0)$, there exists a time $T>0$, and functions $(\rhoteta^{\eps})_{\theta}$, $(\xiteta^{\eps})_{\theta}$ in $C\left([0,T],B_{\delta_1}\right)$, solutions to \eqref{multiphasic} with initial conditions $(\rhoteta(0))_{\theta}$, $(\xiteta(0))_{\theta}$, $(E^0_\eps, B^0_\eps)$.

Moreover the solutions enjoy the following uniform estimates. There exists $C>0$ such that for all $\eps \in (0,\eps_0]$, 
\begin{align}
\label{analytic-uniform-bounds}
&\underset{\theta}{\sup} |\rhoteta^\eps |_{L^\infty(0,T; B_{\delta_1})} + \underset{\theta}{\sup} |\xiteta^\eps |_{L^\infty(0,T; B_{\delta_1})} \leq C \\
\label{analytic-uniform-bounds-EB}
& |E|_{L^\infty(0,T; B_{\delta_1})}+  |B|_{L^\infty(0,T; B_{\delta_1})} \leq C \eps^{-\gamma}.
\end{align}

\end{Thm}

The sequence of solutions that we have obtained thanks to Theorem~\ref{thm2} are so regular in $x$ (see  Lemma~\ref{lembdelta} below) that all requirements of Definition~\ref{defsolweakvm} are of course satisfied.
The uniform initial controls required in Theorem~\ref{theo1} also follow: \eqref{growth-fields} with $\gamma_1= \gamma_2 = \gamma$ is a consequence of~\eqref{analytic-uniform-bounds-EB} (and of Lemma~\ref{lembdelta}).
We also have,  by~\eqref{analytic-uniform-bounds},
$$
\ba
\left\| \rhovm_\eps \right\|_{L^\infty(0,T; L^\infty(\T^3))} &= \left\| \int_{M}\rhoteta^\eps(t,x) \mu(d\theta)  \right\|_{L^\infty(0,T; L^\infty(\T^3))}  \\
&\leq \underset{\theta}{\sup}  \| \rhoteta^\eps \|_{L^\infty(0,T_0; L^\infty(\T^3))}\leq C,
\ea
$$
and likewise, for $\alpha=1$,
$$
\ba
 &\left\|\int_{\bR^3}|v(\xi)| \fvm_\eps(t,x,d\xi)\, \right\|_{L^{\infty} (0,T;L^1(\T^3))} \leq \left\|\int_{\bR^3}|\xi| \fvm_\eps(t,x,d\xi)\, \right\|_{L^{\infty} (0,T;L^\infty(\T^3))}  \\
&\qquad \qquad\qquad= \left\| \int_{M}\rhoteta^\eps(t,x) |\xiteta^\eps(t,x)| \mu(d\theta)  \right\|_{L^\infty(0,T; L^\infty(\T^3))}  \\
&\qquad \qquad\qquad \leq \underset{\theta}{\sup}  \| \rhoteta^\eps \|_{L^\infty(0,T; L^\infty(\T^3))} \underset{\theta}{\sup}  \| \xiteta^\eps \|_{L^\infty(0,T; L^\infty(\T^3))} \leq C,
\ea
$$
that corresponds to $\beta=0$.
\\

Section~\ref{sec-prooftheo2} is dedicated to the proof of Theorem~\ref{thm2}. In order to be able to apply Theorem~\ref{theo1}, we need to check that the associated Vlasov-Poisson solution is suitable, which is done in Theorem~\ref{thm-analytic-vp} in Section~\ref{sec-analytic-vp}.

\section{Proof  of Theorem~\ref{thm2}}
\label{sec-prooftheo2}

In this section we prove an existence result to the multifluid system $\eqref{multiphasic}$ by adaptating the proof of Grenier in the paragraph $2$ of \cite{Grenier-96}, which is itself an adaptation of a simplified proof of the Cauchy-Kovalevskaya Theorem due to Caflisch \cite{Caflisch-90}.

\subsection{Definitions, notations and preliminary results}
To set up our Cauchy-Kovalevskaya argument we first need to consider a scale of Banach spaces $B_{\delta}^{\eta}$.

\begin{Def}
Let $0< \beta < 1$, $\delta_0>1$ and $\eta>0$ be fixed, we consider the following Banach space
$$
B_{\delta_0}^{\eta}:= \left\{u \in C^0([0, \eta(\delta_0 -1)]\times\bT^3), \, \forall \, 0 \leq t \leq \eta(\delta_0-1), \quad u(t) \in B_{\delta_0 - \frac{t}{\eta}}\right\},
$$
endowed with the norm:
$$
\|u\|_{\delta_0}=: \underset{\substack{1 \leq \delta \leq \delta_0, \\ 0 \leq t \leq \eta(\delta_0 - \delta)}}{\sup}\left(|u(t)|_{\delta} + \left(\delta_0 - \delta - \frac{t}{\eta}\right)^\beta |\nabla_x u(t)|_{\delta}\right).
$$
\end{Def}

The space $\bdeltazero$ is a space of functions that are continuous with respect to time with values into the set of analytic functions over the torus $\bT^3$ which takes into account loss of analyticity (in other words, the shrinking of the analyticity domain) as time goes by. Time is as a result bounded by the parameter $\eta(\delta_0-1)$. In the following we are going to prove a local existence result thanks to an iteration scheme, and will consider the parameter $\eta$ as a small parameter.\\

We list some lemmas that are useful for the analysis, whose proofs except the very last one are detailed in Section $2.2$ of \cite{Grenier-96}. The proof of the last lemma is postponed to the Appendix.

\begin{Lem}\label{lembdelta}
For all $\delta,\delta' >1$,
\begin{itemize}
    \item $B_{\delta} \subset B_{\delta'}$ if $\delta' \leq \delta$, 
    \item $\forall s \in \bR$, $B_{\delta}\subset H^s $, the map being compact.
    \item $B_{\delta}$ is  a Banach algebra. Moreover, for $f, g \in B_\delta$, we have
    $$
    |fg|_\delta \leq |f|_\delta |g|_\delta.
    $$
\end{itemize}
\end{Lem}

The advantage of the norms $|.|_{\delta}$ lies in particular in the following lemma: loosely speaking, the $|.|_{\delta'}$ norm of the derivative of a function $f$ in $B_{\delta}$, $1 < \delta' < \delta$, can be controlled by the $|.|_{\delta}$ norm of $f$.

\begin{Lem}
Let $\delta>1$, if $f\in B_{\delta}$ then for any $1<\delta' < \delta$, $i\in\{1,2,3\}$,
$$
|\d_i f |_{\delta'} \leq \dfrac{\delta}{\delta - \delta'} |f|_{\delta}.
$$
\end{Lem}
We have related results around the  space $\bdeltazero$.

\begin{Lem}\label{lembdeltazero}
\begin{itemize}
    \item If $f$ and $g$ are in $\bdeltazero$, then $fg$ is in $\bdeltazero$ as well and 
    $$
    \|fg\|_{\delta_0} \leq \|f\|_{\delta_0}\|g\|_{\delta_0},
    $$
    and in particular if  $\delta + \frac{t}{\eta} < \delta_0$,
    $$
    |\nabla(fg)|_{\delta} \leq \|f\|_{\delta_0}\|g\|_{\delta_0} \left(\delta_0-\delta -\frac{t}{\eta} \right)^{-\beta}.
    $$
   \item if $f$ is in $\bdeltazero$, and if $\delta + \frac{t}{\eta} < \delta_0$,
   \begin{equation*}
| \partial^2_{i,j} f(t)|_\delta \leq C \|f\|_{\delta_0} \delta_0 \left(\delta_0 - \delta - \frac{t}{\eta}\right)^{-\beta-1}.
\end{equation*}

    \end{itemize}
\end{Lem}
   Finally, we have
\begin{Lem}\label{analytic}
If $h$ is analytic and can be written as
$$
h(z) = \sum_{n=0}^{+\infty} a_n z^n,
$$
for $z$ in the disk of center $0$ and radius $R$, $B(0,R)$, and if $f$ is in $\bdelta$ with $|f|_{\delta}< R$, then $h(f)$ is in $\bdelta$, and 
$$
|h(f)|_{\delta} \leq \sum_{n=0}^{+\infty} |a_n| |f|_{\delta}^n.
$$
\end{Lem}

\subsection{Estimate on the force field $F_L$}
Ultimately, we are going to set up an iterative scheme, therefore we will need some a priori analytic bounds on the different quantities that show up in the equation. We begin with the Lorentz force.

Before starting, let us state a useful consequence of Lemma~\ref{analytic}.
 
\begin{Lem}
\label{lem-vxi}
If 
$$\underset{0 \leq t \leq \eta(\delta_0 - \delta)}\sup \underset{\theta\in M }\sup |\xiteta |_{\delta_0} \leq \frac{1}{\sqrt{2}\eps}$$ 
then there is $C>0$ such that for all $\theta \in M$,
\begin{equation} \label{vxiteta}
 \|v(\xiteta)\|_{\delta_0} \leq C \|\xiteta\|_{\delta_0}. %\leq 2 |\xiteta|_{\delta_0} .
\end{equation}

\end{Lem}

\begin{proof} We have  for $\delta \in (1,\delta_0]$, by Lemma~\ref{lembdelta}, 
\begin{align*}
|v(\xi_\theta)|_\delta&= \left|\sum_{n=0}^{+\infty} \dfrac{1/2\times (1/2-1)\times...\times(1/2-(n-1))}{n!}\eps^{2n}|\xi_\theta|^{2n} \xi_\theta \right|_\delta \\
&\leq \sum_{n=0}^{+\infty} \left|\dfrac{1/2\times (1/2-1)\times...\times(1/2-(n-1))}{n!}\right|\eps^{2n}|\xi_\theta|_\delta^{2n+1} \\
&\leq \dfrac{|\xi_\theta|_\delta}{\left(1- \eps^2 |\xi_\theta|_\delta^2\right)^{1/2}} \\
&\leq \sqrt{2} |\xi_\theta|_\delta,
\end{align*}
where we have used the bound on $(\xi_\theta)$ to conclude.

The other part of the estimate likewise follows, according to the formula
$$
\partial_x v(\xi_\theta)=  \frac{\partial_x \xi_\theta}{(1+\eps^2 |\xi_\theta|^2)^{1/2}}  
+  \frac{\eps^2 \xi_\theta (\xi_\theta \cdot \partial_x \xi_\theta)}{(1+\eps^2 |\xi_\theta|^2)^{3/2}}.
$$
\end{proof}

In the following we shall accordingly  systematically assume the following uniform bound:
$$
\underset{0 \leq t \leq \eta(\delta_0 - \delta)}\sup \underset{\theta\in M }\sup |\xiteta |_{\delta_0} \leq \frac{1}{\sqrt{2}\eps}. 
$$

We recall that each phase labeled by $\theta$ is accelerated by a Lorentz force $F_L(t,x,\xiteta)$ produced by the electromagnetic field $(E,B)$, produced collectively by all the phases, and we have
$$
F_{L,\theta} (t,x,\xiteta):= -\nabla_x \phi(t,x) - \eps \d_t A(t,x) +\eps v(\xiteta)\times \left(\nabla\times A(t,x) + \langle B^0 \rangle\right),
$$
with
\begin{equation}\label{poissonteta}
-\Delta_x \phi = \int_{M}\rhoteta \mu(d\theta) -1 ,
\end{equation}
and
\begin{equation}\label{wafeAteta}
\eps^2\partial_t^2 A-\Delta_x A=\eps \bP \left(\int_{ M}v(\xiteta(t,x))\rhoteta(t,x) \mu(d\theta)\right).
\end{equation}
We introduce the quantity
$$
G_{L,\theta}(t,x,\xiteta) = \int_0^t F_L(s,x,\xiteta) ds,
$$
because to obtain a priori estimates on the phase density $\rhoteta$ and the phase momentum field $\xiteta$, if we have a closer look at the equations \eqref{multiphasic}, one can see that we need to control the force field integrated with respect to  time. This is reminiscent of the characteristic equations \eqref{charvm} for the Vlasov-Maxwell system where the time derivative of the velocity of the particles are driven by the Lorentz force $F_L$.

\begin{Lem}
\label{lem-GL}
Assume $\underset{0 \leq t \leq \eta(\delta_0 - \delta)}\sup \underset{\theta\in M }\sup |\xiteta |_{\delta_0} \leq \frac{1}{\sqrt{2}\eps}$.
There exists a positive constant $C$ depending only on the parameters $\beta$ and $\delta_0$ (and not on $\eps$ nor $\eta$) such that
\be\label{gtetabound}
\ba
\|G_{L,\theta}\|_{\delta_0} &\leq   C  \eta  \, \underset{\theta}{\sup} \|\rhoteta -1 \|_{\delta_0}\\
&+ C \eta \left(1+\underset{\theta}{\sup}\|\xiteta\|_{\delta_0}\right) \left(\underset{\theta}{\sup}\|\xiteta\|_{\delta_0}  \underset{\theta}{\sup} \|\rhoteta\|_{\delta_0}\right)  \\
&+ C \eps \left(1+ \underset{\theta}{\sup}\|\xiteta\|_{\delta_0}\right)\left(\|B^0\|_{\delta_0}+ \| E^0\|_{\delta_0} +  \underset{\theta}{\sup} \| \rho_\theta(0)\|_{\delta_0}\right).\ea
\ee
and such that if we consider two solutions $(\rhoteta^,,\xiteta^1)$ and $(\rhoteta^2,\xiteta^2)$ to $\eqref{multiphasic}$, we also have the following stability estimate
\be\label{gtetastab}
\ba
&\|G_{L,\theta}^1 - G_{L,\theta}^2\|_{\delta_0} \leq  C \eta  \,  \underset{\theta}{\sup} \|\rhoteta^1 - \rhoteta
^2 \|_{\delta_0} \\
& + C \eta  \left(1+ \left(\underset{\theta}{\sup}\|\xiteta^1\|_{\delta_0}+ \underset{\theta}{\sup}\|\xiteta^2\|_{\delta_0}\right)\right) \bigg[\left(\underset{\theta}{\sup}\|\xiteta^1\|_{\delta_0}+ \underset{\theta}{\sup}\|\xiteta^2\|_{\delta_0}\right)\underset{\theta}{\sup}\|\rhoteta^1-\rhoteta^2\|_{\delta_0} \\
& +\left(\underset{\theta}{\sup}\|\rhoteta^1\|_{\delta_0}+ \underset{\theta}{\sup} \|\rhoteta^2\|_{\delta_0}\right)\underset{\theta}{\sup}\|\xiteta^1 - \xiteta^2\|_{\delta_0}\bigg]\\
& +  C \eps \underset{\theta}{\sup}\|\xiteta^1 - \xiteta^2\|_{\delta_0}\left(\|B^0\|_{\delta_0} + \| E^0\|_{\delta_0} +  \underset{\theta}{\sup} \| \rho_\theta(0)\|_{\delta_0}\right).
\ea
\ee
\end{Lem}

\subsubsection{Estimates on $\|\int_0^t \nabla \phi\|_{\delta_0}$.}
First from \eqref{poissonteta}, we have for any $k$ in $\bZ^3\setminus\{0\}$:
$$
\left|\cF(  \phi)(k)\right| =  \frac{1}{|k|^2} \cF\left(\int_{M}\rhoteta \mu(d\theta) -1 \right)(k).
$$
%One can remark that this expression is well defined even for $k=0$ because 
%$$
%\cF\left(\int_{M}\rhoteta \mu(d\theta) -1 \right)(0)=0.
%$$
Lemma \ref{lembdeltazero} implies then that for any $t\leq \eta (\delta_0 -1)$,
\begin{equation}\label{nablaphi}
\left|\int_0^t \nabla_x \phi \right|_{\delta} \leq \int_0^t \dfrac{1}{\left(\delta_0 - \delta - \dfrac{s}{\eta}\right)^{\beta}}ds \left| \int_{M}\rhoteta \mu(d\theta) -1 \right|_{\delta}.
\end{equation}
Likewise, for any $j$ in $\{1,2,3\}$, for any $t \in [0, \eta(\delta_0 -1))$
\begin{equation}\label{djnablaphi}
\left|\int_0^t \d_j\nabla_x  \phi \right|_{\delta} \leq  \int_0^t \dfrac{1}{\left(\delta_0 - \delta - \dfrac{s}{\eta}\right)^{\beta+1}}ds \left| \int_{M}\rhoteta \mu(d\theta) -1 \right|_{\delta}.
\end{equation}
Moreover one has
\begin{equation}\label{beta}
\int_0^t \dfrac{ds}{\left(\delta_0 - \delta - \dfrac{s}{\eta}\right)^{\beta}}= \eta \left[\frac{1}{1-\beta}\left(\delta_0-\delta-\dfrac{s}{\eta}\right)^{1-\beta}\right]^{t/\eta}_0 \leq \dfrac{2\eta}{1-\beta}\delta_0^{1-\beta},
\end{equation}
and
\begin{equation}\label{beta+1}
\int_0^t \dfrac{ds}{\left(\delta_0 - \delta - \dfrac{s}{\eta}\right)^{\beta+1}} = \eta \left[\frac{-1}{\beta}\left(\delta_0-\delta-\dfrac{s}{\eta}\right)^{-\beta}\right]^{t/\eta}_0 \leq \dfrac{2\eta}{\beta}\left(\delta_0-\delta-t\right)^{-\beta}.
\end{equation}
%which is a factor tending to infinity as the time $t$ approaches the lifespan of the solution but it is compensated by the factor in the definition of the norm $\|.\|_{\delta_0}$, which has been manufactured precisely for this purpose.\\
Therefore we have proved
\begin{equation}
\label{electrostaticpart}
\left\|\int_0^t \nabla_x \phi\right\|_{\delta_0}\leq C \eta  \, \underset{\theta}{\sup} \|\rhoteta -1 \|_{\delta_0}.
\end{equation}
Likewise, if we consider two solutions $(\rhoteta^1, \xiteta^1, \phi^{1}, A^1)$ and $(\rhoteta^2, \xiteta^2, \phi^{2}, A^2)$ to \eqref{multiphasic}, it comes
\begin{equation}
\label{stab-electrostaticpart}
\left\|\int_0^t\left(\nabla_x \phi^{1} - \nabla_x \phi^{2}\right)
\right\|_{\delta_0} \leq  C\eta \, \underset{\theta}{\sup} \|\rhoteta^1 - \rhoteta
^2 \|_{\delta_0}.
\end{equation}

\subsubsection{Estimates on $\left\|\int_0^t \eps v(\xi_\theta) \times B \right\|_{\delta_0}$.}

Using  Lemmas \ref{lembdelta}, \ref{lembdeltazero} and~\ref{lem-vxi}, we obtain
\begin{align*}
\left|\left(\int_0^t \eps v(\xiteta) \times B   ds\right)\right|_{\delta}&\leq  C \eps \, \underset{\theta}{\sup}\| v(\xiteta)\|_{\delta_0} \| B\|_{\delta_0} \int_0^t \left(\delta_0-\delta-\frac{s}{\eta}\right)^{-\beta}ds   \\
&\leq C\eps \eta  \, \underset{\theta}{\sup}\|\xiteta\|_{\delta_0} \| B\|_{\delta_0}
\end{align*}
Likewise, for $i$ in $\{1,2,3\}$, using  Lemma \ref{lembdeltazero},  the formula \eqref{beta+1}, and proceeding as for the estimation of $\d_j\nabla\phi$ in \eqref{djnablaphi}:
$$ 
\ba
\left|\int_0^t \eps \d_i\left( v(\xiteta) \times B \right) \, ds\right|_{\delta} &\leq C \eps  \,\underset{\theta}{\sup}\|\xiteta\|_{\delta_0}\|B\|_{\delta_0} \int_0^t \left(\delta_0-\delta-\frac{s}{\eta}\right)^{-1-\beta}ds \\
&\leq C \eps \eta  \, \underset{\theta}{\sup}\|\xiteta\|_{\delta_0}\|B\|_{\delta_0} \left(\delta_0 - \delta - \frac{t}{\eta}\right)^{\beta}.
\ea 
$$
In other words, we have
\begin{equation}
\label{estim-vB}
\left\|\left(\int_0^t \eps v(\xiteta) \times B   ds\right)\right\|_{\delta_0}\leq C \eps \eta   \,\underset{\theta}{\sup}\|\xiteta\|_{\delta_0}\|B\|_{\delta_0}.
\end{equation}
The next natural step  consists in estimating $ \|B\|_{\delta_0}$.

\begin{Lem}
\label{lem-estimB}
Assume $\underset{0 \leq t \leq \eta(\delta_0 - \delta)}\sup \underset{\theta\in M }\sup |\xiteta |_{\delta_0} \leq \frac{1}{\sqrt{2}\eps}$.
The following estimate holds:
\begin{equation}
\label{estim-B}
\ba
\| B\|_{\delta_0} &\leq C \eta  \underset{\theta}{\sup}\|\xiteta\|_{\delta_0}\underset{\theta}{\sup} \|\rho_\theta\|_{\delta_0} \\
&+  \| B^0 -\langle B^0\rangle\|_{\delta_0} + |\langle B^0\rangle| + \| E^0\|_{\delta_0} +  \underset{\theta}{\sup} \| \rho_\theta(0)\|_{\delta_0}  + |\langle B^0\rangle|.
\ea
\end{equation}

\end{Lem}

\begin{proof}

We first solve the wave equation~\eqref{wafeAteta} in the Fourier variable $k\in\bZ^3\setminus\{0\}$:
\be
\ba \label{tranfoA}
 \cF(A)(t,k) &= \int_0^t \frac{1}{|k|}\cF\left(\bP\left(\int_{M}v(\xiteta)\rhoteta \mu(d\theta)\right) \right)(s,k)\sin\left(\frac{|k|}{\eps}(t-s)\right)ds \\
 &+ A|_{t=0}(k)\cos\left(\frac{|k|t}{\eps}\right) \\
 &+  \frac{\eps}{|k|}\sin\left(\frac{|k|}{\eps}t\right) \cF\left(\d_t A|_{t=0}\right)(k).
\ea
\ee
Consequently, for any $k\in\bZ^3\setminus\{0\}$,
\be \label{nablaA}
\ba
\cF(\nabla_x\times A)&(t,k) \\
&= \int_0^t \frac{1}{|k|}\cF\left(\nabla_x\times   \bP\left(\int_{M}v(\xiteta)\rhoteta \mu(d\theta)\right) \right)(s,k)\sin\left(\frac{|k|}{\eps}(t-s)\right)ds \\
&+ \nabla_x\times  A^0(k)\cos\left(\frac{|k|t}{\eps}\right) \\
 &+  \frac{\eps}{|k|}\sin\left(\frac{|k|}{\eps}t\right) \cF\left(\nabla_x\times \d_t A|_{t=0}\right)(k).
\ea 
\ee
Let us study the $\| \cdot \|_{\delta_0}$ norm corresponding to these terms. For the first term in~\eqref{nablaA}, we use the fact that (see e.g. the appendix)
$$
\left|\cF(\bP(*))(k)\right|\leq 2 |\cF(*)(k)|.
$$
We can then argue as for the previous estimates to bound its contribution by
$$
C  \eta  \underset{\theta}{\sup}\|\xiteta\|_{\delta_0}  \underset{\theta}{\sup}\|\rhoteta\|_{\delta_0}
$$
The treatment of the contributions of the initial data is straightforward, yielding a bound by
$$
C\left(\| A^0 \|_{\delta_0} + \| \eps \partial_t A|_{t=0}\|_{\delta_0} \right) \leq C  \| B^0 -\langle B^0\rangle\|_{\delta_0}+ \| E^0\|_{\delta_0} + \underset{\theta}{\sup} \| \rho_\theta(0)\|_{\delta_0},
$$
where we have used Lemma~\ref{biotetsavart}. Recalling that
$$
B = \nabla_x \times A + \langle B^0\rangle,
$$
the proof of the lemma is finally complete.

\end{proof}
Gathering~\eqref{estim-vB} and~\eqref{estim-B}, we finally obtain
\be \label{magneticpart}
\ba 
&\left\|\int_0^t \eps v(\xiteta)\times\left(\nabla\times A+ \langle B^0\rangle\right)\right\|_{\delta_0} \leq C \eps \eta \left[\underset{\theta}{\sup}\|\xiteta\|^2_{\delta_0} \, \underset{\theta}{\sup}\|\rhoteta\|_{\delta_0}\right]\\
&\qquad +  C \eps \underset{\theta}{\sup}\|\xiteta\|_{\delta_0}\left(\|B^0\|_{\delta_0}+ \| E^0\|_{\delta_0} +  \underset{\theta}{\sup} \| \rho_\theta(0)\|_{\delta_0} \right).
\ea 
\ee
Similarly, considering two solutions with the same initial data to \eqref{multiphasic} indexed by $i\in\{1,2\}$, we moreover obtain the stability estimate
\be 
\label{stab-magneticpart}
\ba
&\left\|\int_0^t \eps v(\xiteta^1)\times\left(\nabla\times A^1+ \langle B^0\rangle\right) - \eps v(\xiteta^2)\times\left(\nabla\times A^2+ \langle B^0\rangle\right)\right\|_{\delta_0}\\
&\leq  C\eps \eta\Bigg(\left(\underset{\theta}{\sup}\|\xiteta^1\|^2_{\delta_0}+ \underset{\theta}{\sup}\|\xiteta^2\|^2_{\delta_0}\right)\underset{\theta}{\sup}\|\rhoteta^1-\rhoteta^2\|_{\delta_0} \\
&+\left(\underset{\theta}{\sup}\|\rhoteta^1\|_{\delta_0}+\underset{\theta}{\sup}\|\rhoteta^2\|_{\delta_0}\right)\left(\underset{\theta}{\sup}\|\xiteta^1\|_{\delta_0}+ \underset{\theta}{\sup}\|\xiteta^2\|_{\delta_0}\right)\underset{\theta}{\sup}\|\xiteta^1 - \xiteta^2\|_{\delta_0}\Bigg)\\
&+ C \eps \underset{\theta}{\sup}\|\xiteta^1 - \xiteta^2\|_{\delta_0}\left(\|B^0\|_{\delta_0}+ \| E^0\|_{\delta_0} +  \underset{\theta}{\sup} \| \rho_\theta(0)\|_{\delta_0}\right).
\ea 
\ee

\subsubsection{Estimates on $\left\|\int_0^t \eps \d_t A(s,x)ds\right\|_{\delta_0}$.}
As we do not need to estimate directly $\eps \d_t A(s,x)$ but its integral over $(0,t)$,
we actually need to study
$$
\int_0^t \eps \d_t A(s,x)ds = \eps\left(A(t,x) - A^0(x)\right).
$$
We can then use the formula \eqref{tranfoA} to obtain
\be \label{transfodtA}
\ba
\cF\left(A- A^0\right)(t,k) = &\int_0^t \frac{1}{|k|}\cF\left(\bP\left(\int_{M}v(\xiteta)\rhoteta \mu(d\theta)\right) \right)(s,k)\sin\left(\frac{|k|}{\eps}(t-s)\right)ds \\
 &+ A^0(k)\left(\cos\left(\frac{|k|t}{\eps}\right)-1\right) \\
 &+  \frac{\eps}{|k|}\sin\left(\frac{|k|}{\eps}t\right) \cF\left(\d_t A|_{t=0}\right)(k).
\ea 
\ee
We study this term exactly as in the proof of Lemma~\ref{lem-estimB}.
It follows that
\begin{equation}
\label{longitudinalE}
\ba
    \left\|\int_0^t \d_t A(s,x)ds \right\|_{\delta_0} &\leq C \eps \eta \left(\underset{\theta}{\sup}\|\xiteta\|_{\delta_0}\|\rhoteta\|_{\delta_0}\right) \\
    &+  C \eps \left(\|B^0- \langle B^0\rangle\|_{\delta_0} + \| E^0\|_{\delta_0} +  \underset{\theta}{\sup} \| \rho_\theta(0)\|_{\delta_0} \right).
    \ea
\end{equation}
Again if we consider two solutions with the same initial data to \eqref{multiphasic}, we obtain
\be
\label{stab-longitudinalE}
\ba
\left\|\int_0^t  \eps\d_t A^1 \, ds - \int_0^t\eps \d_t A^2 \, ds \right\|_{\delta_0} \leq  & C \eta \Big[\left(\underset{\theta}{\sup}\|\xiteta^1\|_{\delta_0}+ \underset{\theta}{\sup}\|\xiteta^2\|_{\delta_0}\right)\underset{\theta}{\sup}\|\rhoteta^1-\rhoteta^2\|_{\delta_0} \\
&+\left(\underset{\theta}{\sup}\|\rhoteta^1\|_{\delta_0}+ \underset{\theta}{\sup} \|\rhoteta^2\|_{\delta_0}\right)\underset{\theta}{\sup}\|\xiteta^1 - \xiteta^2\|_{\delta_0}\Big].
\ea
\ee

\begin{Rmk}
The gain of a power of $\eps$ in~\eqref{longitudinalE} and~\eqref{stab-longitudinalE} which is due to the integration in time is somehow reminiscent of the treatment of $I_3$ in Theorem~\ref{theo1} (see Section~\ref{sec-I3}).
\end{Rmk}

Gathering~\eqref{electrostaticpart}, \eqref{magneticpart} and~\eqref{longitudinalE},  we find that
\be
\ba
\|G_{L,\theta}\|_{\delta_0} &\leq  C  \eta  \, \underset{\theta}{\sup} \|\rhoteta -1 \|_{\delta_0}\\
&+ C \eta \left(1+\underset{\theta}{\sup}\|\xiteta\|_{\delta_0}\right) \left(\underset{\theta}{\sup}\|\xiteta\|_{\delta_0}  \underset{\theta}{\sup} \|\rhoteta\|_{\delta_0}\right)  \\
&+ C \eps \left(1+ \underset{\theta}{\sup}\|\xiteta\|_{\delta_0}\right)\left(\|B^0\|_{\delta_0} + \| E^0\|_{\delta_0} +  \underset{\theta}{\sup} \| \rho_\theta(0)\|_{\delta_0}\right).
\ea
\ee
and  by~\eqref{stab-electrostaticpart}, \eqref{stab-magneticpart} and~\eqref{stab-longitudinalE} , we deduce the stability estimate
\be
\ba
&\|G_{L,\theta}^1 - G_{L,\theta}^2\|_{\delta_0} \leq  C \eta  \,  \underset{\theta}{\sup} \|\rhoteta^1 - \rhoteta
^2 \|_{\delta_0} \\
& + C \eta  \left(1+ \left(\underset{\theta}{\sup}\|\xiteta^1\|_{\delta_0}+ \underset{\theta}{\sup}\|\xiteta^2\|_{\delta_0}\right)\right) \bigg[\left(\underset{\theta}{\sup}\|\xiteta^1\|_{\delta_0}+ \underset{\theta}{\sup}\|\xiteta^2\|_{\delta_0}\right)\underset{\theta}{\sup}\|\rhoteta^1-\rhoteta^2\|_{\delta_0} \\
& +\left(\underset{\theta}{\sup}\|\rhoteta^1\|_{\delta_0}+ \underset{\theta}{\sup} \|\rhoteta^2\|_{\delta_0}\right)\underset{\theta}{\sup}\|\xiteta^1 - \xiteta^2\|_{\delta_0}\bigg]\\
& +  C \eps \underset{\theta}{\sup}\|\xiteta^1 - \xiteta^2\|_{\delta_0}\left(\|B^0\|_{\delta_0}+ \| E^0\|_{\delta_0} +  \underset{\theta}{\sup} \| \rho_\theta(0)\|_{\delta_0}\right).
\ea
\ee
This concludes the proof of Lemma~\ref{lem-GL}.

\subsection{Estimates for $\rhoteta$ and $\xiteta$}
In this section, we prove some a priori analytic bounds on $\rhoteta$ and $\xiteta$. As in the previous subsections we assume that 
$$
\underset{0 \leq t \leq \eta(\delta_0 - \delta)}\sup \underset{\theta\in M }\sup |\xiteta |_{\delta_0}\leq \frac{1}{\sqrt{2}\eps},
$$
which we recall by Lemma \ref{lem-vxi} implies that
$$
\underset{\theta}{\sup}\|v(\xiteta)\|_{\delta_0} \leq C \underset{\theta}{\sup} \|\xiteta\|_{\delta_0}.
$$
We consider $\widetilde{\xiteta}$ the solution to
$$
\d_t \widetilde{\xiteta} + \left(v(\xiteta)\cdot\nabla_x\right)\xiteta = F_{L,\theta}, \qquad  \widetilde{\xiteta} (0) = \xiteta(0).
$$
By Lemma \ref{lembdelta}, Lemma \ref{lembdeltazero} 
and Lemma~\ref{lem-vxi}, it comes that
\begin{align*}
    |\widetilde{\xiteta}(t)|_{\delta} &\leq \int_0^t |\d_t\widetilde{\xiteta}(s)|_{\delta} ds + |\widetilde{\xiteta}(0)|_{\delta}\\
    &\leq \|{\xiteta}(0)\|_{\delta_0} + C \int_0^t \left(\delta_0-\delta-\frac{s}{\eta}\right)^{-\beta}\|\xiteta\|^2_{\delta_0}ds +  \|G_{L,\theta}\|_{\delta_0} \\
    &\leq \|{\xiteta}(0)\|_{\delta_0} + C \eta  \|\xiteta\|^2_{\delta_0} +\|G_{L,\theta}\|_{\delta_0}.
\end{align*}
The same argument holds for the estimate bearing on $|\d_i\widetilde{\xiteta}(t)|_{\delta}$ and it comes
\begin{equation}\label{xitetabound}
\|\widetilde{\xiteta}\|_{\delta_0} \leq \|\widetilde{\xiteta}(0)\|_{\delta_0} +  C \eta \|\xiteta\|^2_{\delta_0} +\|G_{L,\theta}\|_{\delta_0} .
\end{equation}
Similarly, considering two solutions to \eqref{multiphasic}, this analysis provides the stability estimate
\be \label{xitetastabi}
\ba
\|\widetilde{\xiteta}^1 - \widetilde{\xiteta}^2\|_{\delta_0} \leq & C \eta \left(\|\xiteta^1\|_{\delta_0}+\|\xiteta^2\|_{\delta_0}\right)   \|\xiteta^1 - \xiteta^2\|_{\delta_0} + \|G_{L,\theta}^1- G_{L,\theta}^2\|_{\delta_0}.
\ea
\ee
\\
Now working on the equation for $\rhoteta$, we consider the solution $\widetilde{\rhoteta}$ to the equation
$$
\d_t \widetilde\rhoteta + \nabla_x \cdot(v(\xiteta)\rhoteta) = 0, \qquad  \widetilde{\rhoteta} (0) = \rhoteta(0).
$$
By a similar argument it turns out that
\begin{equation}\label{rhotetabound}
\|\widetilde{\rhoteta}\|_{\delta_0} \leq \|{\rhoteta}(0)\|_{\delta_0} +  C \eta  \|\xiteta\|_{\delta_0}\|\rhoteta\|_{\delta_0},
\end{equation}
and
\be\label{rhotetastabi}
\ba
\|\widetilde{\rhoteta}^1 - \widetilde{\rhoteta}^2\|_{\delta_0} &\leq   C \eta \left((\|\xiteta^1\|_{\delta_0}+\|\xiteta^2\|_{\delta_0}\right)\|\rhoteta^1 - \rhoteta^1\|_{\delta_0} \\
&+  C \eta  (\|\rhoteta^1\|_{\delta_0}+\|\rhoteta^2\|_{\delta_0})\|\xiteta^1-\xiteta^2\|_{\delta_0}.
\ea 
\ee

\subsection{The iterative scheme}
We define inductively $(\rhoteta^n)_{n\in \bN} $, $(\xiteta^n)_{n\in \bN}$, $(G_{L,\theta}^n)_{n\in \bN}$ as follows.

For $n=0$  we set
$$
F_{L,\theta}^0 = 0,
$$
$$
\xiteta^0(t) = \xiteta(0),
$$
and
$$
\rhoteta^0 (t) = \rhoteta(0),
$$
for all $ 0<t < \eta$; for $n\geq 1$,  we rely on the induction relation
\begin{align*}
 &\d_t \rhoteta^{n+1} + \nabla_x \cdot(v(\xiteta^n)\rhoteta^n) = 0,\\
 & \d_t {\xiteta}^{n+1} + \left(v(\xiteta^n)\cdot\nabla_x\right)\xiteta^n ) =F_{L,\theta}^n,\\
 & \text{with  }  \rhoteta^{n+1}(0)=\rhoteta^n(0) \quad \text{ and }\quad \xiteta^{n+1}(0)=\xiteta^n(0),
\end{align*}
in which
$$
F_{L,\theta}^n (t,x,\xiteta):= -\nabla_x \phi^n(t,x) - \eps \d_t A^n(t,x) +\eps v(\xiteta^n)\times \left(\nabla\times A^n(t,x) + \langle B^0 \rangle\right),
$$
where $\phi^n$ and $A^n$ solve
\begin{equation*}
-\Delta_x \phi^n = \int_{M}\rhoteta^n \mu(d\theta) -1 ,
\end{equation*}
and
\begin{equation*}
\eps^2\partial_t^2 A^n-\Delta_x A^n=\eps \bP \left(\int_{\theta \in M}v(\xiteta^n(t,x))\rhoteta^n(t,x) \mu(d\theta)\right),
\end{equation*}
with the initial condition
$$
\nabla \times A^n|_{t=0}= B^0 - \langle B^0\rangle, \quad \eps \partial_t A^n|_{t=0} = E^0 + \nabla_x \phi^0.
$$

\subsubsection{Estimates on $\rhoteta^0$, $\xiteta^0$.}
%From the assumptions made in Theorem~\ref{thm2} on $\xiteta^0$,$\rhoteta^0$, we deduce that $\|\xiteta(0)\|_{\delta_0}$, $\|\rhoteta(0)\|_{\delta_0}$ are bounded uniformly on $\theta$, $\eps$ by a constant $\widetilde{C_1}$ which only depends on the initial conditions.\\

\subsubsection{Contraction estimates for $n\geq 1$}

\begin{Lem}
\label{lem-induction}
There exist $\eta>0$, $C_1,C_2>0$ and $\eps_0>0$  such that if $\eps \in (0,\eps_0]$,
then:
\begin{itemize}
    \item for all $n\geq 0$:
          \be \label{inductionbound}
          \ba 
          &\underset{\theta}{\sup}\|\rhoteta^n\|_{\delta_0}\leq  C_1, \\ %2\|\rhoteta^{0,n}\|_{\delta_0}\leq 2 \tilde{C_1},\\
          &\underset{\theta}{\sup} \|\xiteta^n \|_{\delta_0}\leq 2 C_1, \\% \ 2 \tilde{C_1},\\
         &\underset{\theta}{\sup}\|G_{L,\theta}^n \|_{\delta_0}\leq C_1; % \frac{\tilde{C_1}}{2}.
         \ea
          \ee
          
     \item   for all $n\geq 1$:  
          \be \label{inductionstabi}
          \ba 
          &\underset{\theta}{\sup}\|\rhoteta^n - \rhoteta^{n-1}\|_{\delta_0}\leq \frac{C_2}{2^n},\\
          &\underset{\theta}{\sup} \|\xiteta^n- \xiteta^{n-1}\|_{\delta_0}\leq \frac{C_2}{2^n},\\
         &\underset{\theta}{\sup} \|G_{L,\theta}^n - G_{L,\theta}^{n-1}\|_{\delta_0}\leq \frac{C_2}{2^{n+2}} .
         \ea
          \ee
\end{itemize}
\end{Lem}

\begin{proof}

Let us first focus on the first item, that is~\ref{inductionbound}. We argue by induction. For $n=0$, we can choose 
$$C_1:= 4C_0,$$ 
which enforces~\eqref{inductionbound} according of the assumptions on the initial data.  We then pick $\eps_0>0$ such that $\eps_0 \leq \dfrac{1}{\sqrt{2}{C_1}}$, which will enable us to apply Lemma~\ref{lem-vxi} and thus all the estimates of the previous subsections are valid.

Let assume that~\eqref{inductionbound} holds for some $n \geq 0$. Then  using  \eqref{gtetabound} and \eqref{inductionbound}, it follows that
\be
\ba
\|G_{L,\theta}^{n+1}\|_{\delta_0} &\leq  C  \eta  \, \underset{\theta}{\sup} \|\rhoteta^n -1 \|_{\delta_0}\\
&+ C \eta   \left(1+\underset{\theta}{\sup}\|\xiteta^n\|_{\delta_0}\right) \left(\underset{\theta}{\sup}\|\xiteta^n\|_{\delta_0}  \underset{\theta}{\sup} \|\rhoteta^n\|_{\delta_0}\right)  \\
&+ C \eps \left(1+ \underset{\theta}{\sup}\|\xiteta^n\|_{\delta_0}\right)\left(\|B^0\|_{\delta_0}  + \| E^0\|_{\delta_0} +  \underset{\theta}{\sup} \| \rho_\theta(0)\|_{\delta_0}\right) \\
&\leq  \eta  C \left( C_1 + 1 + (2C_1+1)2C_1^2 \right) + 3\eps_0^{1-\gamma}  C (1+2C_1)  C_0.
\ea
\ee
Recalling that $\gamma \in [0,1)$, choosing $\eta$ and $\eps_0$ sufficiently small, we get
$$
\|G_{L,\theta}^{n+1}\|_{\delta_0} \leq C_1.
$$
Similarly, using \eqref{inductionbound} and the estimates \eqref{xitetabound}, we  obtain
\be 
\ba 
\|\xiteta^{n+1}\|_{\delta_0} \leq & \|\xiteta^0\|_{\delta_0} + C \eta \|\xiteta^n\|_{\delta_0}^2 + \|G_{L,\theta}^n\|_{\delta_0}\\
\leq & C_0 + \eta C C_1^2 + C_1 \\
\leq & 2 C_1,
\ea 
\ee
up to taking $\eta$ small enough.
We omit the treatment of $\|\rhoteta^{n+1}\|_{\delta_0}$ which is completely similar.  We have therefore proved by induction \eqref{inductionbound}.

We now prove~\eqref{inductionstabi} by induction. The case $n=1$ requires a special treatment. We actually use the rough bounds
$$
\ba
&\underset{\theta}{\sup}\|\rhoteta^1 - \rhoteta^{0}\|_{\delta_0}\leq \underset{\theta}{\sup}\|\rhoteta^1\|_{\delta_0}  +   \underset{\theta}{\sup} \| \rhoteta^{0}\|_{\delta_0} \leq   {2C_1}, \\
&\underset{\theta}{\sup}\|\xiteta^1- \xiteta^{0}\|_{\delta_0}\leq \underset{\theta}{\sup}  \|\xiteta^1\|_{\delta_0}  +\underset{\theta}{\sup}  \| \xiteta^{0}\|_{\delta_0} \leq 4 C_1, \\
&\underset{\theta}{\sup} \|G_{L,\theta}^1 - G_{L,\theta}^{0}\|_{\delta_0}\leq \underset{\theta}{\sup} \|G_{L,\theta}^1 \|_{\delta_0}\leq C_1,
\ea
 $$
 and we choose $C_2 := 2^3 C_1$, so that~\eqref{inductionstabi} holds for $n=1$.
 
Assume now~\eqref{inductionstabi} holds for some $n\geq 1$.
Using \eqref{inductionbound} and \eqref{gtetastab} we obtain
\be
\ba
&\|G_{L,\theta}^{n+1} - G_{L,\theta}^{n}\|_{\delta_0} \leq  C \eta  \,  \underset{\theta}{\sup} \|\rhoteta^{n} - \rhoteta
^{n-1} \|_{\delta_0}  + C \eta  \left(1+ \left(\underset{\theta}{\sup}\|\xiteta^{n}\|_{\delta_0}+ \underset{\theta}{\sup}\|\xiteta^{n-1}\|_{\delta_0}\right)\right) \\
&\qquad\times\bigg[\left(\underset{\theta}{\sup}\|\xiteta^{n}\|_{\delta_0}+ \underset{\theta}{\sup}\|\xiteta^{n-1}\|_{\delta_0}\right)\underset{\theta}{\sup}\|\rhoteta^{n}-\rhoteta^{n-1}\|_{\delta_0} \\
&\qquad \qquad\qquad\qquad +\left(\underset{\theta}{\sup}\|\rhoteta^{n}\|_{\delta_0}+ \underset{\theta}{\sup} \|\rhoteta^{n-1}\|_{\delta_0}\right)\underset{\theta}{\sup}\|\xiteta^{n} - \xiteta^{n-1}\|_{\delta_0}\bigg]\\
&\quad +  C \eps \underset{\theta}{\sup}\|\xiteta^{n} - \xiteta^{n-1}\|_{\delta_0}\left(\|B^0\|_{\delta_0}+ \| E^0\|_{\delta_0} +  \underset{\theta}{\sup} \| \rho_\theta(0)\|_{\delta_0}\right) \\
&\leq  \eta C ( 1 + (1+  2C_1)2C_1 + 2C_1)\frac{C_2}{2^n} + 3\eps_0^{1-\gamma}C C_0   \frac{C_2}{2^n}\\
&\leq \frac{C_2}{2^{n+3}},
\ea
\ee
up to taking $\eta$  and $\eps_0$ small enough.
Likewise, using \eqref{inductionbound}, \eqref{inductionstabi} and \eqref{xitetastabi} we obtain
\begin{align*}
    \|\xiteta^{n+1}- \xiteta^n\|_{\delta_0} &\leq C  \eta \left(\|\xiteta^{n}\|_{\delta_0}+\|\xiteta^{n-1}\|_{\delta_0}\right)   \|\xiteta^{n} - \xiteta^{n-1}\|_{\delta_0} + \|G_{L,\theta}^{n}- G_{L,\theta}^{n-1}\|_{\delta_0}\\
   &\leq \left(4 \eta  C C_1 + \frac{1}{2} \right)    \frac{C_2}{2^{n}} \\
    &\leq \dfrac{C_2}{2^{n+1}},
\end{align*}
up to taking $\eta$ small enough.
We argue similarly for $\|\rhoteta^{n+1} - \rhoteta^n\|_{\delta_0}$, which allows to close the induction argument.

\end{proof}

To conclude, Lemma~\ref{lem-induction} proves that for all $\theta \in M$, $(\rhoteta^n)_n$, $(\xiteta^n)_n$ are Cauchy sequences in the Banach spaces $\bdeltazero$ for a suitable small parameter $\eta$. As a result, they  converge to functions $\rhoteta$ and $\xiteta$ belonging to $\bdeltazero$. Letting $n$ tend to infinity one can check that the pair $(\rhoteta, \xiteta)_\theta$ is a solution to the system \eqref{multiphasic}. 
Now let $\delta_1 \in (1,\delta_0)$. We  pick $T = \eta (\delta_0 -\delta_1)$ to conclude the existence part of the theorem.
\\

There remains to derive the claimed uniform in $\eps$ estimates.
The uniform bound in $\eps$ for $(\rhoteta)_\theta$, $(\xiteta)_\theta$ is a consequence of~\eqref{inductionbound}. The control of $(E,B)$ is a consequence of the formula~\eqref{tranfoA}, of the formula for $\eps \partial_t A$ 
\be
\ba \label{tranfodtA}
 &\cF(\eps \partial_t A)(t,k) =\\
  &\int_0^t \frac{\eps}{|k|}\cF\left(\bP\left(\int_{M}\nabla_\xi v(\xiteta) :\partial_t \xiteta \rhoteta + v(\xiteta) \partial_t \rhoteta\mu(d\theta)\right) \right)(s,k)\sin\left(\frac{|k|}{\eps}(t-s)\right)ds \\
 &+\frac{1}{|k|} \sin\left(\frac{|k|}{\eps}t\right)\cF\left(\bP\left(\int_{M} v(\xiteta^0) \rhoteta^0\mu(d\theta)\right) \right)(k)\\
 &+ \cF(A|_{t=0})(k)|k| \sin\left(\frac{|k|t}{\eps}\right) \\
 &+  \cos\left(\frac{|k|}{\eps}t\right) \cF(\eps \partial_t A|_{t=0})(k),
\ea
\ee
of the assumption on the initial electromagnetic field~\eqref{E0B0init} and of the above uniform bounds for $(\rhoteta)_\theta$, $(\xiteta)_\theta$.

\subsection{The Vlasov-Poisson case}
\label{sec-analytic-vp}

In order to be able to apply Theorem~\ref{theo1}, we also need to check that there exists a suitable weak solution to the Vlasov-Poisson system (in the sense of Definition~\ref{sol-vp}) associated to the initial condition
$$
f^0(x, d\xi) = \int_{M}\rhoteta(0,x)\delta(\xi - \xiteta(0,x)) \mu(d\theta).
$$
We recall that according to \cite{Loeper-2006}, such a suitable weak solution is then unique.
As in the Vlasov-Maxwell case, we look for the solution under the form
$$
\fvp(t,x,\xi) = \int_{M}\rhoteta(t,x)\delta(\xi - \xiteta(t,x)) \mu(d\theta),
$$
with $(\rho^\theta, \xi^\theta)$ solving the multifluid system
\be\label{multiphasic-vp}
\left\{
\ba
 &\forall \theta\in M,\, \d_t \rhoteta(t,x) +\nabla_x \cdot\left( \xiteta(t,x) \rhoteta(t,x)\right)=0  \,,
\\
&\forall \theta\in M ,\, \d_t\xiteta (t,x) + \left(\xiteta(t,x)\cdot\nabla_x\right)\xiteta(t,x)= -\nabla_x \phi, \\
 &- \Delta \phi = \int_{M}\rhoteta(t,x) \mu(d\theta) - 1 \,,\\
&\forall \theta\in M, \quad \rhoteta|_{t=0}=\rhoteta(0), \quad \xiteta|_{t=0}=\xiteta(0).
\ea
\right.
\ee
We can obtain the following result.
\begin{Thm}
\label{thm-analytic-vp}
Let $(M,\mu)$ the probability space used to define the multifluid system \eqref{multiphasic-vp}, let $\delta_0> \delta_1>1$, let $C_0>0$, and $(\rhoteta(0))_{\theta}$, $(\xiteta(0))_{\theta}$ be  families of $\bdeltazero$ such that
\begin{equation}
\begin{aligned}
&\underset{\theta}{\sup} |\rhoteta (0) |_{\delta_0} \leq C_0, \\
&\underset{\theta}{\sup} |\xiteta(0) |_{\delta_0} \leq C_0.
\end{aligned}
\end{equation}
Then there exists a time $T_0>0$, and functions $(\rhoteta)_{\theta}$, $(\xiteta)_{\theta}$ in $C\left([0,T_0],B_{\delta_1}\right)$, solutions to \eqref{multiphasic-vp} with initial conditions $(\rhoteta(0),\xiteta(0))_{\theta}$.
\end{Thm}
We will not give the proof of this result as it is already contained in that of Theorem~\ref{thm2} (see also \cite{Grenier-96}).
We can check that this solution is suitable. Indeed, 
$$
\ba
\left\| \rhoVp \right\|_{L^\infty(0,T_0; L^\infty(\T^3))} &= \left\| \int_{M}\rhoteta(t,x) \mu(d\theta)  \right\|_{L^\infty(0,T_0; L^\infty(\T^3))}  \\
&\leq \underset{\theta}{\sup}  \| \rhoteta \|_{L^\infty(0,T_0; L^\infty(\T^3))}<+\infty,
\ea
$$
and likewise
$$
\ba
 &\left\|\int_{\bR^3}|\xi|^4\fvp(t,x,d\xi)\, \right\|_{L^{\infty} (0,T_0;L^1(\T^3))}  \\
&\qquad \qquad\qquad= \left\| \int_{M}\rhoteta(t,x) |\xiteta(t,x)|^4\mu(d\theta)  \right\|_{L^\infty(0,T_0; L^1(\T^3))}  \\
&\qquad \qquad\qquad \leq \underset{\theta}{\sup}  \| \rhoteta \|_{L^\infty(0,T_0; L^\infty(\T^3))} \underset{\theta}{\sup}  \| \xiteta \|^4_{L^\infty(0,T_0; L^\infty(\T^3))}<+\infty.
\ea
$$
This finally shows that we can apply Theorem~\ref{theo1} to such solutions (note though that to completely enter the framework of Theorem~\ref{theo1}, one also needs to enforce all conditions on the initial data of Definition~\ref{cond-initial}).

%\section{Conclusion and remarks}

\section{Appendix}
\subsection{On the Leray projection $\bP$}
We gather in this paragraph some remarks on the Leray projection $\bP: L^2(\bT^3)\rightarrow L^2(\bT^3)$, its Fourier transform and its continuity properties with respect to any norm $H^s$, $s\in [0,+\infty)$.\\

More generally let us explain how a vector field can be decomposed in a divergence free part and an irrotational part (the Helmholtz decomposition). Given $F\in C^1_c (\bT^3, \bT^3)$ one can define $\bP(F)$
 as follows:
 $$
 \bP(F)(x) = F(x) + \nabla \psi(x), \quad \text{ for any}\, x \in \bT^3,
 $$
 with $\psi$ satisfying the equation
 $$
 -\Delta \psi = \nabla \cdot F.
 $$
 One can check that it implies that
 $$
 \nabla_x\cdot \bP(F) =0.
 $$
 Applying the Fourier transform, one obtains for any $k\in\bZ^3$:
 $$
 \cF(\bP(F))(k)= \cF(F)(k)+ ik \cF(\psi),
 $$
 and 
 $$
 \cF(\psi)(k) = \dfrac{i}{|k|^2}k\cdot\cF(F)(k),
 $$
 which implies 
 \begin{equation}\label{bP}
 \cF(\bP(F))(k) = \left(Id-\dfrac{k\otimes k}{|k|^2}\right)\cF(F)(k).
\end{equation}
 Therefore 
 $$
 |\cF(\bP(F))(k)| \leq C  |\cF(F)(k)|,
 $$
 and by the Plancherel Theorem $\bP$ extends to a continuous operator on $L^2(\bT^3)$ characterized by the formula $\eqref{bP}$. From the same formula we have the continuity with respect to any Sobolev norm $H^s$, $s\in[0,+\infty)$.
 
 \subsection{Proof of Lemma \ref{analytic}}
 We recall that all the proofs concerning the properties of the analytic norms we use can be found in Section $2.2$ of \cite{Grenier-96}. Let us  however explain Lemma \ref{analytic}.\\
 
One can define the function $h(f)$ over $\bT^3$  by the power series:
$$
h(f)(x)= \sum_{n=0}^{+\infty}a_n f(x)^n.
$$
Then we have the following inequality, thanks to Lemma \ref{lembdelta} according to which the space $\bdelta$ is a Banach algebra,
$$
|h(f)|_{\delta} \leq \sum_{n=0}^{+\infty}|a_n|\left|f^n\right|_{\delta},
 $$
 and we can conclude.

\bibliographystyle{plain}
\bibliography{myrefs}

\end{document}